\documentclass{amsart}
\usepackage{amsfonts, amsmath, amstext, amssymb,hyperref}

\usepackage{color,centernot}

\usepackage[dvips]{graphicx}
\usepackage{latexsym}
\newtheorem{thm}{Theorem}[section]
\newtheorem{prop}[thm]{Proposition}

\newtheorem{lemma}[thm]{Lemma}

\newtheorem{cor}[thm]{Corollary}

\newtheorem{question}[thm]{Question}

\newtheorem{definitiontemp}[thm]{Definition}
\newenvironment{defn}{\begin{definitiontemp}
\normalfont}{\end{definitiontemp}}

\def\bec{\begin{cor}}
\def\enc{\end{cor}}

\def\bet{\begin{thm}}
\def\ent{\end{thm}}

\def\becor{\begin{cor}}
\def\encor{\end{cor}}

\def\bel{\begin{lem}}
\def\enl{\end{lem}}

\def\bedef{\begin{defn}}
\def\endef{\end{defn}}

\def\bep{\begin{prop}}
\def\enp{\end{prop}}

\newenvironment{pf}{\begin{trivlist}\item[\hskip\labelsep
{\it Proof.}]}{\end{trivlist}}

\newcommand{\set}[2]{\ensuremath{ \{ #1 : #2 \} }}

\renewcommand{\deg}[1]{\ensuremath{\text{deg}(#1)}}

\newcommand{\N}{\mathbb{N}}
\newcommand{\Z}{\mathbb{Z}}
\newcommand{\hatZ}{\widehat{\mathbb{Z}}}
\newcommand{\F}{\mathfrak{F}}

\newcommand{\Q}{\mathbb{Q}}
\newcommand{\R}{\mathbb{R}}

\renewcommand{\S}{\mathcal{S}}
\newcommand{\C}{\mathbb{C}}
\newcommand{\A}{\mathcal{A}}
\newcommand{\B}{\mathcal{B}}

\renewcommand{\L}{\mathcal{L}}

\newcommand{\avec}{\vec{a}}

\newcommand{\cvec}{\vec{c}}
\newcommand{\fvec}{\vec{f}}
\newcommand{\Fvec}{\vec{F}}
\newcommand{\Gvec}{\vec{G}}

\newcommand{\pvec}{\vec{p}}
\newcommand{\Pvec}{\vec{P}}

\newcommand{\x}{\text{$\bf{x}$}}

\renewcommand{\c}{\text{$\bf{c}$}}
\newcommand{\xvec}{\vec{x}}

\newcommand{\Xvec}{\vec{X}}
\newcommand{\yvec}{\vec{y}}
\newcommand{\Yvec}{\vec{Y}}

\newcommand{\cross}{^{\textsf{x}}}

\newcommand{\la}{\langle}
\newcommand{\ra}{\rangle}

\def\converges{\!\downarrow}

\newcommand{\dom}[1]{\text{dom}(#1)}

\newcommand{\bfd}{\boldsymbol{d}}
\def\bfz{\boldsymbol{0}}
\def\s01{\ensuremath{\Sigma^0_1}}
\def\d02{\ensuremath{\Delta^0_2}}
\def\phi{\varphi}

\def\res{\!\!\upharpoonright\!\!}

\newcommand{\comment}[1]{}

\begin{document}

\title{Computability for tree presentations of continuum-size structures}

\author{Jason Block}
\address{
\parbox{\dimexpr\textwidth -\parindent}{ \hspace*{0.1in}  Department of Mathematics, College of William \& Mary, 200 Ukrop Way, Williamsburg, VA 23185 USA}}
\email{Jeblock@wm.edu}
\urladdr{\url{https://sites.google.com/view/jasonblockmath/}}
\author{Russell Miller}
\address{
\parbox{\dimexpr\textwidth -\parindent}{ \hspace*{0.1in}  Department of Mathematics, Queens College -- City University of New York, 65-30 Kissena Blvd., Flushing, NY 11367 USA; and \\
\hspace*{0.1in}  Ph.D. Programs in Mathematics and Computer Science, Graduate Center - City University of New York, 365 Fifth Avenue, New York, NY 10016 USA}}
\email{Russell.Miller@qc.cuny.edu}
%\urladdr{\url{https://qcpages.qc.cuny.edu/~rmiller/}}
%}

\begin{abstract}
We formalize an existing computability-theoretic method of presenting
first-order structures whose domains have the cardinality of the continuum.
Work using these methods until now has emphasized their topological properties.
We shift the focus to first-order properties, using computable structure theory
(on countable structures) as a guide.  We present three basic questions
to be asked when a structure is presented as the set of paths through a
computable tree, as in our definition, and also propose the concept of
\emph{tree-decidability} as an analogue to the notion of decidability
for a countable structure.  As examples, we prove decidability results
for certain additive and multiplicative groups of $p$-adic integers, products of these
(such as the profinite completion of $\Z$), and the field of real numbers.
\end{abstract}

\maketitle

\section{Tree presentations of structures}
\label{sec:treepres}

Computable structure theory has long focused on countable structures.
In the standard notion of a \emph{computable structure} $\mathcal C$,
domain elements are represented by natural numbers (thus requiring
the structure to be countable), so that all $n$-ary functions and relations
in $\mathcal C$ map $\omega^n$ into either $\omega$
(for functions) or $\{0,1\}$ (for relations).  One then requires that
these maps all be computable, in the usual Turing sense of functions on $\omega$.

Continuum-sized structures clearly require a more expansive definition.
In Definition \ref{defn:treepres} we will propose one: a \emph{tree presentation}
of a structure, along with variants that will be introduced subsequently.
This concept is based largely on existing notions from computable analysis \cite{W00},
and has already been used by several research groups studying continuum-sized
structures, often with a distinct flavor of computable topology (see the book
\cite{DM25} for a comprehensive treatment; also \cite{DM23},\cite{K23}, \cite{P20} and \cite{W09}).  
While topology is indeed essential in this context,
we hope to add the perspective of computable structure theory
and enable mathematicians to pose, about continuum-sized structures,
the sorts of model-theoretic questions that have traditionally been asked about
computable (countable) structures.  In this section, after defining tree
presentations, we will list several such questions.

To introduce tree presentations, we recall the additive group $\Z_p^+$ of the
\emph{$p$-adic integers}.  Two equivalent definitions exist, for each prime $p$.
In the first definition, a $p$-adic integer $x$ is simply an $\omega$-tuple $(j_n)_{n\geq 0}$
from $(\Z/p\Z)^\omega$, usually represented as
$$ x = \sum_{n\geq 0} j_n\cdot p^{n}.$$
Here the group addition is essentially base-$p$ addition.
In the second definition, a $p$-adic integer $x$ is defined to be an $\omega$-tuple
$$ x = (k_1,k_2,k_3,\ldots) \in (\Z/p\Z) \times (\Z/p^2\Z)\times (\Z/p^3\Z)\times\cdots$$
satisfying $(\forall n) k_{n+1}\equiv k_n\bmod p^n$,
with the group operation $+$ being coordinatewise addition
(modulo $p^n$ in the $n$-th coordinate).
This second definition will correspond to the ``tree presentation'' $T_p^+$
of $\Z_p^+$ that we describe later in this section.

In each of these definitions, the elements of $\Z_p^+$ are given naturally as paths through
the complete $p$-branching tree $T_p$.  The first definition
names $p$-adic integers simply as paths in the tree $(\Z/p\Z)^{<\omega}$.
Moreover, the group operation
is effective as a function from $[T_p]^2$ into $[T_p]$,
being computable by a Turing functional.  (See \cite[III.1]{S87} or \cite[2.4.1]{DH10}
for this concept.)  The identity element is just the
constant path $0^\omega$.  The second definition is similarly effective,
although the nodes in the tree are superficially different.
Before going into it, we need to make the basic concept rigorous.

From these examples and others, we draw our concept of a
\emph{tree presentation} of a structure, as defined here.  Recall that $[T]$
denotes the set of (infinite) paths through a subtree $T$ of $\omega^{<\omega}$,
and that a node in $T$ is \emph{extendible} if it lies on a path in $[T]$.
\begin{defn}
\label{defn:treepres}
Fix an oracle $C\subseteq\omega$, and let $\L$ be a
$C$-computable signature with no relation symbols except equality.
A \emph{$C$-computable tree presentation} of an $\L$-structure consists of:
\begin{itemize}
\item
a $C$-computable subtree $T\subseteq \omega^{<\omega}$
in which all nodes are extendible; and
%\item
%for each constant symbol $c$ in $\L$, a $D$-computable path $p_c\in [T]$;
\item
for each $n$-ary function symbol $f$ in $\L$, a Turing functional $\Phi_f$
such that, for each $f$ and every $(X_1,\ldots,X_n)\in[T]^n$,
the function $\Phi_f^{C\oplus X_1\oplus\cdots\oplus X_n}:\omega\to\omega$ is total and lies in $[T]$.
\end{itemize}

Constant symbols $c$ are treated as $0$-ary function symbols,
with $\Phi_c^C:\omega\to\omega$ defining a path through $T$.
In an infinite signature $\L$, we also require that the indices of the functionals $\Phi_f$
be $C$-computable uniformly in $f$.

We write $\A_T$ for the $\L$-structure with domain $[T]$ thus presented:
$$ \A_T = ([T]; \Phi_f^{C\oplus\cdots}~|~f\in\L).$$
If $\S \cong \A_T$, we call $\A_T$ a \emph{($C$-computable) tree presentation of $\S$}.
\end{defn}

Definition \ref{defn:treepres} can be strengthened or weakened in various ways.
The weaker notion of a \emph{nonextendible tree presentation} employs the same
definition, except that $T$ is allowed to contain nonextendible nodes.  In the present work
we will not meet any such presentations, but they have appeared, for example,
in the work \cite{L81} of La Roche on computable Galois groups.  Indeed, in \cite{S81}
Smith constructed a computable nonextendible tree presentation of a profinite group
that has no computable extendible tree presentation.  On the other hand,
one can strengthen Definition \ref{defn:treepres} significantly by requiring $T$
to be finite-branching (as holds in both \cite{L81} and \cite{M24}) and/or by requiring
that the branching in $T$ be $C$-computable (so that $C$ computes the function
$\sigma\mapsto |\set{\tau\in T}{|\tau|=|\sigma|+1~\&~\sigma\sqsubset\tau}|$).

The tree presentation $T_p^+$ of $\Z_p^+$ that we began to describe above
satisfies all of these stronger requirements, with $C=\emptyset$.
$T_p^+$ is the following tree:
$$\set{(k_1,\ldots,k_n)\in\omega^{<\omega}}{(\forall i)0\leq k_i<p^i~\&~
(\forall i <n) k_{i+1}\equiv k_i\bmod p^i},$$
with $(j_i)_{i>0}+(k_i)_{i>0} = (j_i+k_i\bmod p^i)_{i>0}$.
One should notice that the identity element is a computable path in $T_p^+$,
and that the unary operation of inversion in the group is also computable.
(In traditional computable structure theory
this is trivial, as one can search through the domain $\omega$
for the inverse.  With uncountable domains, such search procedures
no longer apply in general.  The reader may wish to prove that
for computably-finite-branching tree presentations such as $T_p^+$,
computability of inversion does follow from computability of the group operation.)

In a traditional computable presentation $\B$ (of a countable structure), the entire
atomic diagram of $\B$ is computable.  In a computable tree presentation,
the functions on $\A_T$ are all computable (now by Turing functionals),
but equality is generally not computable this way.  This is simply a fact of life
when dealing with uncountable structures.
%% Likewise, when we define \emph{decidable tree presentations} in Section \ref{sec:decidable}, equality will have to be given rather than decided.
\begin{lemma}
\label{lemma:equality}
In a tree presentation $T$, if there exists an oracle $C$
such that, for all $P,Q\in [T]$, $\Phi^{C\oplus P\oplus Q}(0)$ halts and outputs
$1$ if $P=Q$ and $0$ otherwise, then $\A_T$ is countable.

Conversely, every $C$-computable countable structure $\B$ has a $C$-computable tree presentation
in which equality is $C$-computable as above.
\end{lemma}
\begin{pf}
Computations are finite, so when $\Phi^{C\oplus P\oplus P}(0)\converges=1$,
some $\sigma\sqsubseteq P$ isolates $P$ in $[T]$.  As $T$ is countable, so must $[T]$ be.

For the converse, use the subtree $\set{(n,0,0,0,\ldots,0)}{n\in\dom{\B}}$ of $\omega^{<\omega}$,
with the path $n\widehat{~}0^{\infty}$ representing the element $n$ of $\B$.
\qed\end{pf}

Definition \ref{defn:treepres} should be compared with \cite[Defn.\ 2.4.4]{DM25}
by Downey and Melnikov,
which is based on their definition of a presentation of a Polish space.
The approach there is grounded in topology, which makes those definitions
appear substantially different from ours.  In fact, though, they are quite closely related.
Topology plays a far greater role in structures of size continuum than it
does for countable structures in traditional computable structure theory,
and so it is quite reasonable for a comprehensive
study of computability on uncountable structures to begin that way.
Here we do not aim to present such a study, and so we are content
with the more model-theoretic Definition \ref{defn:treepres} for the present.
Nevertheless, it is certainly wise, when building a tree presentation of $\A$,
to ensure that it reflects the topology on $\A$,
whenever there is a natural topological structure to reflect.
In particular, the topology on $\A_T$ should coincide with
the usual topology on its domain $[T]$, given by the basis of sets
$$ \set{P\in [T]}{\sigma\sqsubset P}$$
as $\sigma$ ranges through all the nodes of $T$.

This topology is totally disconnected, which makes it futile to attempt to give topologically
appropriate tree presentations of structures such as the field $\R$ of real numbers.
Nevertheless, the usual method of presenting $\R$ in computable analysis,
using fast-converging Cauchy sequences from $\Q$, does fit into the framework
of Definition \ref{defn:treepres} with just a small modification.  Consider the following
subtree of $\Q^{<\omega}$:
$$ T_{\R} = \left\{(q_1,q_2,\ldots,q_n)\in\Q^{<\omega}~:~(\forall k\leq n)(\forall j<k)~|q_j-q_k|\leq\frac1{2^j}\right\}.$$
A path in $T_{\R}$ is a \emph{fast-converging Cauchy sequence} of rational numbers,
for the condition $|q_j-q_k|\leq\frac1{2^j}$ gives convergence and implies the condition 
$(\forall j)~|q_j-\lim_sq_s|\leq \frac1{2^j}$ that defines \emph{fast convergence}.
Moreover, every $r\in\R$ is the limit of some path in $[T_\R]$.
The only reason why this is not a tree presentation of $\R$ is that each single
real number is the limit of many different paths.  The next definition is a generalization
accommodating this issue.

\begin{defn}
\label{defn:treequotientpres}
Fix an oracle $C\subseteq\omega$, and let $\S$ be a structure in a
$C$-computable signature $\L$ with no relation symbols except equality.
A \emph{$C$-computable tree quotient presentation} of $\S$ consists of:
\begin{itemize}
\item
a $C$-computable subtree $T\subseteq \omega^{<\omega}$
in which all nodes are extendible;
%\item
%for each constant symbol $c$ in $\L$, a $D$-computable path $p_c\in [T]$;
\item
for each function symbol $f$ in $\L$, a Turing functional $\Phi_f$; and
\item
an equivalence relation $\sim$ on $[T]$ that is $\Pi^0_1$ 
relative to $C$,
\end{itemize}
such that, for each ($n$-ary) $f$ and every $(X_1,\ldots,X_n)\in[T]^n$,
the function $\Phi_f^{C\oplus X_1\oplus\cdots\oplus X_n}:\omega\to\omega$ is total and lies in $[T]$,
and such that the quotient structure
$$ \A_T=([T]/\!\sim~; \Phi_f^{C\oplus\cdots})$$
is isomorphic to $\S$.
In particular, every function $\Phi_f^{C\oplus\cdots}$ must respect $\sim$,
so as to yield a well-defined map from $([T]/\!\sim)^n$ into $[T]/\!\sim$.
\end{defn}

It is quickly seen that the $T_\R$ given above is a computable tree quotient presentation
of the field $\R$, as the operations of addition and multiplication on fast-converging
Cauchy sequences have long been known to be computable by Turing functionals.
(So likewise are subtraction and division, apart from the issue of dividing by zero.)
The equivalence relation $\sim$ on two sequences $(q_n)$ and $(r_n)$
says that they have the same limit, and because of their fast convergence
this can be expressed as $(\forall n)~|q_n-r_n|\leq\frac2{2^n}$, which is $\Pi^0_1$.

We find tree presentations to be intuitively more clear than tree quotient presentations.
In the same sense, in model theory, an interpretation of one first-order structure $\A$
in another one $\B$ (given by a definable subset of some power $\B^n$,
a definable equivalence relation on this subset, and formulas defining the
relations and functions of $\A$ there) is more pleasing if the equivalence relation
is pure equality, in which case it is called a definition of $\A$ in $\B$,
rather than merely an interpretation.
%is not as pleasing as a definition (i.e., an interpretation using
%equality as the equivalence relation on the domain of the copy of $\A$ inside $\B$).
Nevertheless, for the purposes of computation, tree quotient presentations are
just as effective as tree presentations.  The reason is the bugaboo of tree presentations:
equality of paths is itself only a $\Pi^0_1$ property, not a decidable property.  Indeed, two paths
$X,Y\in[T]$ satisfy $X=Y$ just if $(\forall n)~X(n)=Y(n)$, and it is quickly seen that
this is in general the best one can do.

For reasons related to these, we have avoided relation symbols (except $=$)
in our signatures.  The natural definition would require that, in a $C$-computable
tree presentation of $\A$, some Turing functional $\Phi_R$, given oracle
$(C\oplus P_1\oplus\cdots\oplus P_n)$, would decide membership of the tuple $\vec{P}$
in the $n$-ary relation $R$, by outputting either $0$ or $1$.  However, only
nearly-trivial relations can be decidable this way, particularly in the context
of most of this article, in which the trees in question are finite-branching:
consequently $[T]$ is compact, and so there must be a finite level $l$ in $T$
such that the membership of $\vec{P}$ in $R$ is determined by the restrictions
$(P_1\res l,\ldots,P_n\res l)$.  Section \ref{sec:decidable} might point towards
possible resolutions of this issue.  A different resolution, not to be considered
in this article, involves continuous logic; see \cite{H23} for an introduction to that topic.

Along with the examples $T_{\R}$ and $T_p^+$ above, another computable tree presentation
of a familiar continuum-sized structure appears in \cite{M24}, in which one of us
considers computability for the absolute Galois group of the rational numbers.
In \cite{B24,B25}, the other of us then carried out similar investigations for profinite
groups more generally.
Those articles propose three basic questions about computable tree presentations $T$
and the corresponding structures $\A_T$, all of which we will examine here for $\Z_p^+$.

\begin{question}
\label{question:elem}
For each Turing ideal $I$, the set $(\A_T)_I=\set{X\in[T]}{\deg{X}\in I}$ forms a substructure
of $\A_T$.  To what extent are such $(\A_T)_I$ elementary substructures of $\A_T$,
or at least elementarily equivalent to $\A_T$?  Of special interest
for profinite groups are the cases where $I=\{\boldsymbol{0}\}$  and where $I$ is a Scott ideal.
\end{question}

Using $T_p^+$,
we will show in Section \ref{sec:elementary} that, for every Turing ideal $I$,
$(\Z_p^+)_I$ is an elementary subgroup of $\Z_p^+$.  (In contrast, \cite{B25}
builds a computably-tree-presentable group where this fails.)

When $T$ is $T_{\R}$, it is known that the fields $\R_I$ are all elementary subfields of $\R$,
as they are all real closed and Th$(\R)$ is model-complete.  In \cite{K99},
Kapoulas proved a similar result for the field $\Q_p$, the $p$-adic completion of $\Q$,
given by Cauchy sequences from $\Q$ that converge fast in the $p$-adic topology.
Further questions are possible:  for example, Korovina and Kudinov showed
in \cite{KK17} that the Turing degree spectrum of the (countable) field $\R_{\{\bfz\}}$ is
precisely the set $\set{\bfd}{\bfz''\leq_T\bfd'}$ of all high degrees.

\begin{question}
\label{question:definableset}
Each formula $\alpha(x_1,\ldots,x_n)$ defines a subset $S_\alpha\subseteq [T]^n$.
What is the computational complexity of such a subset?
\end{question}
For a $\Sigma_n$ formula $\alpha$, the quantification is over elements of $[T]$,
so this appears to be a question about levels in the analytic hierarchy:
the complexity of $S_\alpha$ is at worst $\Sigma^1_n$.  However,
in many cases it turns out to be much lower than this, often arithmetical.
Indeed, we do not yet know of any example of such a set $S_\alpha$
that lies outside the arithmetical hierarchy.  For $\Z_p^+$ and also
for $\R$, the complexity is far lower, as described in Section \ref{sec:decidable}.
This holds even when $\alpha$ is allowed to use parameters and one
considers the complexity of $S_\alpha$ uniformly in the parameters.

\begin{question}
\label{question:Skolem}
Given a structure $\A$ and a $\Sigma_n$ formula $\alpha(x_1,\ldots,x_n)$ of the form 
$\exists  y \beta(\vec x, y)$,
a \emph{(generalized) Skolem function for $\alpha$}
is a partial function $F:\A^n\to \A$ such that:
\begin{itemize}
    \item the domain of $F$ includes all tuples $\vec a$ such that $\A\models \alpha(\vec a)$; and
    \item if $\A\models \alpha (\vec a)$, then $\A\models \beta(\vec a, F(\vec a))$.
\end{itemize}
Thus, such an $F$
picks out a witness for $\alpha(\vec a)$ whenever one exists. Given a $\vec b$ such that $\A\models \neg \alpha(\vec b)$,  either $\vec b$ will not be in the domain of $F$ or we will simply have  $\A\models \neg \beta(\vec b,F(\vec b))$. 
(If $\A\models\forall\xvec~\exists y~\alpha(\xvec, y)$, then this $F$
is a Skolem function in the usual sense.)

For a given tree presentation $T$ and formula $\alpha$, a Skolem function $F$ for $\alpha$ in $T$ is computable if there exists a Turing functional $\Psi$ such that for all $\vec P$ in the domain of $F$,  we have that $\Psi^{\vec P}:\omega\to \omega$ is total and equal to the path $F(\vec P)\in [T]$.  Our question  is whether
a computable Skolem function exists for a given $\alpha$; or more generally, how low the complexity
of a Skolem function for $\alpha$ can be.  Computing such an $F$ might require
the jump $(\Pvec)'$ as input, for example, or several jumps over $\Pvec$.
One can also ask to what extent Skolem functions for different formulas $\alpha$
can be computed uniformly in $\alpha$.
\end{question}

For purely existential formulas $\exists y~\beta(\xvec,y)$, Question \ref{question:Skolem} can be seen
as asking whether the situation for computable (countable) structures remains true.
In a computable structure, for arbitrary parameters $\pvec$ from the domain $\omega$,
one can simply search for a  $y$ from $\omega$ satisfying $\beta(\pvec,y)$:
the search will find such a tuple if any exists and will then output the $y$ it found.
In a tree-presented structure, the analogous situation
would be that all existential formulas have computable Skolem functions uniformly in the formula.
By Proposition \ref{prop:or issue}, this is never completely true, but in Section
\ref{sec:existential} we will see that it can hold to a significant extent.

For Question \ref{question:definableset}
 above, it is important to have a working definition of
complexity for subsets of $[T]$.  Full definitions appear in (e.g.)\ \cite{K95},
but for our purposes, it suffices to consider
$S\subseteq[T]^k$ to be $\Sigma^0_{n+1}$ just if there exists a Turing functional
$\Phi$ such that
$$ (\forall (X_1,\ldots,X_k)\in[T]^k)~[(X_1,\ldots,X_k)\in S \iff \Phi^{(X_1\oplus\cdots\oplus X_k)^{(n)}}(0)\converges]$$
and to be $\Pi^0_{n+1}$ just if $([T]^k-S)$ is $\Sigma^0_{n+1}$.
(Here $A^{(n)}$ denotes the $n$-th Turing jump of a set $A\subseteq\omega$.)
For a fixed oracle $C$, the condition for $S$ to be $(\Sigma^0_{n+1})^C$
%(abbreviated $\Sigma^S_k$ if arithmeticity is clear)
is that $\Phi^{(C\oplus X_1\oplus\cdots\oplus X_k)^{(n)}}(0)$
should halt for precisely the elements of $S$;
and $S$ is (boldface) $\bf{\Sigma^0_{n+1}}$ if it is $(\Sigma^0_{n+1})^C$
for some $C\subseteq\omega$.
Definition \ref{defn:effectivelyopen} will introduce the term
\emph{effectively open} for (lightface) $\Sigma^0_1$ sets.

\section{Elementarity of subgroups}
\label{sec:elementary}

For $\Z_p^+$,
we first address Question \ref{question:elem}, about the elementarity of
subgroups, which turns out to have a fairly short answer:
every Turing ideal defines an elementary subgroup of $\Z_p^+$.  In contrast,
\cite{B25} builds a computable tree presentation of a profinite group
with the opposite answer, and for the absolute Galois group of $\Q$
it is open whether the computable elements form an elementary
(or even elementarily equivalent) subgroup.

The proof of our Proposition \ref{prop:full elemen} will follow from results
of Szmielew in \cite{S55}, formulated more recently by Eklof in \cite{E72}.
Let $G$ be an additive group. We say that $G$ is $\emph{of finite exponent}$ if every element of $G$ has finite order. Given $n\in \N$ define $nG=\{ng:g\in G\}$ and $G[n]=\{g\in G:ng=0\}$. 

\begin{thm}[\cite{S55}, see Theorem 1 in \cite{E72}]
\label{thm:E72 1}
If $A$ and $B$ are abelian groups, then $A\equiv B$ if and only if
\begin{align}
    A \text{ is of finite exponent } \iff B \text{ is of finite exponent;}
\end{align}
and for all primes $q$ and all $n\in \N$:
\begin{align}
    \dim (q^n A[q]/q^{n+1}A[q])&= \dim (q^n B[q]/q^{n+1}B[q]),\\
    \lim_{n\to \infty}\dim(q^n A/q^{n+1}A)&=  \lim_{n\to \infty}\dim(q^n B/q^{n+1}B),\\
    \lim_{n\to \infty}\dim(q^n A[q])&=  \lim_{n\to \infty}\dim(q^n B[q]),
\end{align}
  where $\dim$ means dimension over $\Z/q\Z$.   
\end{thm}

\begin{thm}[\cite{S55}, see Theorem 2 in \cite{E72}]
\label{thm:E72 2}
Given an abelian group $B$, a subgroup $A$ is an elementary subgroup  of $B$ if and only if $A\equiv B$ and $nA=nB\cap A$ for all $n\in \N$.
    
\end{thm}

\begin{prop}
\label{prop:full elemen}
    For every Turing ideal $I$, the subset
$$ (\Z_p^+)_I = \set{g\in\Z_p}{\deg{g}\in I}$$
forms an elementary subgroup of $\Z_p^+$.
\end{prop}

\begin{pf}
It is immediate that $(\Z_p^+)_I$ is a subgroup, as the sum of any two elements
of this subset is computable (indeed uniformly) from the join of the degrees
of the two elements, while the additive inverse of $g$ is Turing-equivalent to $g$ itself. Additionally, the identity element of $\Z_p^+$ is computable and the degree $\bold{0}$ is in every ideal $I$. 
For the claimed elementarity, it is clear that $n (\Z_p^+)_I= n \Z_p^+ \cap (\Z_p^+)_I $,
so we need only show that $(\Z_p^+)_I\equiv \Z_p^+$. 

No non-identity element of $\Z_p^+$ has finite order, so condition (1) of Theorem \ref{thm:E72 1} is met. Given any prime $q$, $(\Z_p^+)_I[q]$ and $\Z_p^+[q]$ are both trivial, so conditions (2) and (4) are both met. If $q\neq p$, then $q^m (\Z_p^+)_I=(\Z_p^+)_I$ and $q^m \Z_p^+=\Z_p^+$ for all $m\in \N$. Thus,  $$\lim_{n\to \infty}\dim(q^n (\Z_p^+)_I/q^{n+1}(\Z_p^+)_I)=0=  \lim_{n\to \infty}\dim(q^n \Z_p^+/q^{n+1}\Z_p^+)$$ for all primes $q\neq p$. Last, note that $$p^n (\Z_p^+)_I/p^{n+1}(\Z_p^+)_I\cong \Z/p\Z \cong p^n \Z_p^+/p^{n+1}\Z_p^+$$ which gives $$\lim_{n\to \infty}\dim(p^n (\Z_p^+)_I/p^{n+1}(\Z_p^+)_I)=1=  \lim_{n\to \infty}\dim(p^n \Z_p^+/p^{n+1}\Z_p^+)$$ and so condition (3) is met as well. \qed

\end{pf}

\section{Deciding Truth in $\Z_p^+$}
\label{sec:existential}

In this section we focus on the additive group $\Z_p^+$,
presented using the computable tree $T_p^+$ defined above, whose paths
are the elements of $\Z_p$, with Turing functionals for addition and negation
that simply perform these operations on the $n$-th coordinate modulo $p^n$.
%This is a good example of a group that (given its size)
%requires a tree presentation, but is presentable as nicely as possible among
%computably tree-presentable groups.  
One cannot expect the atomic diagram
of $\Z_p^+$ to be decidable, as the mere task of deciding whether two particular
paths are equal is clearly undecidable.
In this section we will show that the difficulty of deciding truth
of more complicated formulas is no worse than that for equality;
and in Section \ref{sec:decidable} we will discuss the difficulty itself.

\begin{defn}
\label{defn:effectivelyopen}
A subset $X\subseteq (\Z_p)^j$ is \emph{effectively open}, or $\Sigma^0_1$,
if it is the union of a c.e.\ collection of basic clopen sets; or equivalently,
if there exists a c.e.\ set $\{(\vec a_k,l_k)\}_{k\in \N} \subseteq \Z^{j+1}$
(an \emph{enumeration of $X$}) such that
$$X=\bigcup_{k\in \N} \{\vec f \in (\Z_p)^j: f_i(l_k)= a_{k,i} \text{ for all }i=1,...,j\}.$$
$X$ is \emph{effectively clopen} if both it and its complement in $(\Z_p)^j$
are effectively open.  By compactness, this is equivalent to $X$ being a finite union as above.
\end{defn}

\begin{defn}
\label{defn:existential-atomic}
A formula of the form $\exists G\alpha(\Fvec,G)$, where $\alpha$
is atomic, is said to be \emph{existential-atomic}.  It is
\emph{properly existential-atomic} if no atomic formula without
$G$ is equivalent (under the axioms for a group) to $\alpha(\Fvec,G)$.
\end{defn}

\begin{lemma}
\label{lemma:Zplinear}
Each single equation in $\Z_p^+$, in variables $F_1,\ldots,F_n,G$,
can be re-expressed in the form $\sum_i a_i\cdot F_i = b\cdot G$
with coefficients $a_i,b$ all in $\Z$.  Assuming $b\neq 0$, the equation
admits a solution just if $p^e$ divides $\sum_i a_i\cdot f_i(e)$, where
$e$ is maximal such that $p^e$ divides $b$.  Thus the set
$\set{\fvec}{(\exists G)~\sum_i a_i\cdot f_i = b\cdot G}$ is effectively clopen,
uniformly in the coefficients $a_i$ and the coefficient $b$ except when $b=0$.
Moreover, there is a procedure that always computes a solution $g$ (provided one exists)
uniformly in the oracle $\fvec$ and also uniformly in the coefficients.  
When $b\neq 0$ that solution $g$ is unique whenever it exists,
with $g\res n$ determined by $\fvec\res n$ whenever $n>e$.
%, uniformly for every tuple $\fvec$ from our tree presentation
%$T_p^+$ of $\Z_p^+$, whether the equation has a solution in $\Z_p^+$.
\end{lemma}
(When $b=0$, the formula $(\exists G) \sum_i a_i\cdot F_i = b\cdot G$
is equivalent to a quantifier-free statement about $\Fvec$.  So Lemma \ref{lemma:Zplinear}
may be viewed as a decision procedure for sets defined by
properly existential-atomic formulas.)
\begin{pf}
Since $\Z_p^+$ is abelian, the re-expression of equations is immediate.
Write $F=\sum a_iF_i$, so the equation is simply $F=bG$.

Of course, if $b=0$, then either there is no solution
(if the given $f=\sum a_if_i$ is nonzero) or else every $g\in\Z_p^+$ is a solution.
So our procedure for computing a solution first checks whether $b=0$
and, if so, simply outputs the path $(0,0,\ldots)$ as $g$.  For this we do not need
to be able to decide whether $f=0$.

Now assume $b\neq 0$.  
Let $e\geq 0$ be maximal with $p^e$ dividing $b$.  For each $n>0$,
let the integer $c_{n}$ be the
multiplicative inverse of $\frac b{p^e}$ in $(\Z/p^{n}\Z)\cross$.  
(Consequently $c_{n+k}\equiv c_n\bmod p^n$ for every $n$ and $k$.)
Now if $p^e$ fails to divide $f(e)$,
then the equation is unsatisfiable:  every $g\in\Z_p$ has
$bg(e)\equiv 0\not\equiv f(e)\bmod p^e$.

Assuming $p^e$ does divide $f(e)$ it also divides $f(n+e)$ since
$f(n+e)\equiv f(e)\equiv 0\bmod p^e$.  So we may define $g\in\Z_p$ by
$$ g(n) = c_{n+e}\cdot\frac {f(n+e)}{p^e}\bmod p^{n}
%~~~(\text{so also~} g(n)= c_{n+e}\cdot\frac {f(n+e)}{p^e}\bmod p^{n+e})
.$$
For each $n$, we can write $f(n+e+1)=f(n+e)+kp^{n+e}$ with $k\in\Z$, so that
\begin{align*}
g(n+1) &\equiv c_{n+e+1}\cdot\frac {f(n+1+e)}{p^e}
\equiv c_{n+e+1}\cdot\frac {f(n+e)+kp^{n+e}}{p^e} \bmod p^n\\
&\equiv c_{n+e}\cdot \frac {f(n+e)}{p^e} +c_{n+e+1}kp^{n}
\equiv g(n) \bmod p^n,
\end{align*}
as required in order for $g$ to lie in $\Z_p$.
%(As inverses of $\frac b{p^e}$ in $\Z/p^{n+e+1}\Z$ and $\Z/p^{n+e}\Z$, $c_{n+e+1}\equiv c_{n+e}\bmod p^n$.)
Thus, for every $m\geq 0$, we have 
$g(e+m)\equiv g(m) \bmod p^m$ and consequently $bg(e+m)\equiv bg(m) \bmod p^{m+e}$
(since $p^e$ divides $b$).  Therefore
\begin{align*}
bg(m+e) &\equiv bg(m)\bmod p^{m+e}\\
&\equiv b c_{m+e} \cdot\frac {f(m+e)}{p^e}\bmod p^{m+e}\text{~~~(again since $p^e|b$)}\\
&\equiv \left(\frac b{p^e}\cdot c_{m+e}\right)\cdot f(m+e)\bmod p^{m+e}\\
&\equiv f(m+e)\bmod p^{m+e}~~~~~~~\text{(by the definition of $c_{m+e}$)}.
\end{align*}
Also, $bg(n)\equiv 0\equiv af(n)\bmod p^n$ for all $n < e$,
because $p^e$ divides both $bg(e)$ and $f(e)$.  Thus $b\cdot g= f$ in $\Z_p^+$,
and the computation of this $g$ was effective, uniformly in the original equation.
(The uniformity includes the case $b=0$, since an algorithm can check whether
$b=0$ and, if so, output any $g$ it likes.)

Finally, notice that when $b\neq 0$ and $p^e$ divides $f(e)$
(hence also divides $f(e+n)$), this solution $g$ is unique,
because the equations above determine every value $g(n)$ within $\Z/p^{n}\Z$.
%(which in turn determine $g(0),\ldots,g(e)$).
Indeed, $bg(e+n) \equiv f(e+n) \bmod p^{e+n}$
so $p^{e+n}$ divides $(bg(e+n)-f(e+n))$. %, which equals $p^e\cdot (\frac b{p^e}g(e+n) - \frac{f(e+n)}{p^e})$.
Thus $p^n$ divides the integer $(\frac b{p^e}g(e+n) - \frac{f(e+n)}{p^e})$,
so $\frac b{p^e}g(e+n) \equiv \frac{f(e+n)}{p^e}\bmod p^n$,
Since $\frac b{p^e}$ is a unit mod $p^n$, this determines $g(e+n)\bmod p^n$,
hence also determines $g(n)\bmod p^n$.
\qed\end{pf}

\begin{cor}
\label{cor:cofinite}
$\Z_p^+\models\forall\Fvec~\exists! G~\sum_i a_i\cdot F_i = b\cdot G$
whenever $p\!\not|~\!b$.
%the prime $p$ does not divide $b$.
\qed\end{cor}

In the language of Question \ref{question:Skolem},
Lemma \ref{lemma:Zplinear} shows that in the tree presentation $T_p^+$,
existential-atomic formulas have computable generalized Skolem functions,
uniformly in the formula (and even when $b=0$, so the $\exists$-quantifier
is superfluous).  When $b\neq 0$, the lemma also
answers Question \ref{question:definableset}.

\begin{defn}
\label{defn:conjunctive}
An existential formula $\beta$ is $\emph{conjunctive}$ if it is of the form
$$\exists \vec G \bigwedge_i \alpha_i(\vec F,\vec G) $$
%\emph{conjunctive existential (with parameters)} if it is of the form
%$$\exists G_1\cdots\exists G_k~\wedge_i \alpha_i(\Fvec,\Gvec),$$
where %$X$ is either $\exists$ or $\forall$ according to the parity of $n$ and 
every $\alpha_i$ is a literal:  either atomic or the negation
of an atomic formula.  We naturally refer to these as \emph{equations}
and \emph{inequations}, respectively.  An equation is \emph{removable}
if no $G_i$ appears in it, and the negation of a removable equation
is a \emph{removable inequation}.  (If $\alpha_1$ is removable, then
$\exists\Gvec(\alpha_0\wedge\alpha_1)$ is equivalent to
$\alpha_1\wedge\exists\Gvec\alpha_0$; hence the name.)

A conjunctive existential sentence
$\exists G_1\dots\exists G_k\alpha(\Fvec,\Gvec)$
is \emph{reduced} if:
\begin{itemize}
\item
its conjuncts include no removable literals, and 
%\item
%each equation in $\alpha$ uses at most one of the variables $G_i$, and
\item
each $G_i$ appears in at most one equation in $\alpha$, and
\item
for each $i$, if $G_i$ appears in an equation in $\alpha$,
then it does not appear in any inequation in $\alpha$.
\end{itemize}
(So $\alpha$ may consist of exactly $k$ equations, with no inequations;
or else of fewer than $k$ equations along with arbitrarily many inequations.)
\end{defn}
\begin{lemma}
\label{lemma:reduced}
Every conjunctive existential formula $\exists\Gvec\alpha(\Fvec, \Gvec)$ is equivalent to 
a conjunction $\beta(\Fvec)\wedge\exists\Gvec\alpha'(\Fvec,\Gvec)$, where $\beta$
is a quantifier-free conjunction and $\exists\Gvec\alpha'$ is reduced.  Moreover, such a
conjunction may be found effectively.
\end{lemma}

\begin{pf}
The key to the proof is that, since the ring $\Z_p$ is an integral domain,
the solutions to an equation $\sum a_if_i =bG$ in $\Z_p^+$
are exactly the solutions to each nonzero integer multiple
$\sum a_icf_i =bcG$ of that equation.  Likewise, the solutions to
$\sum a_if_i \neq bG$ are exactly those to $\sum a_icf_i \neq bcG$.

First, therefore, if some $G_i$ occurs in an equation $\alpha_m$ 
(saying $bG_i=A_m$ for some $\Z$-linear combination $A_m$ of $\Fvec$ and the other $G_j$'s)
and also occurs in some other equation or inequation $\alpha_n$
(saying that $b'G_i$ does or does not equal such a combination $A_n$),
we replace $\alpha_n$ by the (in)equation $\alpha_n'$ given by
$b'A_m=bA_n$ (or $b'A_m\neq bA_n$), in which $G_i$ no longer appears,
but which yields the same solutions as the original system.
Repeat this wherever possible, so that we satisfy the final two conditions for being reduced.
Then move all removable literals outside the $\exists$-quantifiers.
This is the procedure required by the lemma.
\qed\end{pf}

Reduced form is not unique:  two reduced formulas can define the same set.
However, the notion as stated here suffices for our purposes.
When $\exists G\alpha$ is in reduced form we often write $\alpha^=$ for
the conjunction of its equations (if any) and
$\alpha^{\neq}$ for the conjunction of its inequations (if any).

Here we encounter a curiosity of this (or any other) computable tree presentation of $\Z_p^+$.
The proof of Lemma \ref{lemma:Zplinear} shows that the properly existential-atomic diagram
of $\Z_p^+$ is decidable, uniformly in the parameters $\fvec$, but the full existential diagram fails
to be decidable, because the truth of atomic sentences with parameters
(which are existential sentences, in a trivial way) cannot be decided.
The undecidability is quite clear for the formula $F=F+F$, which
holds for exactly one path $f$ (the identity) through our tree presentation of $\Z_p^+$
and fails for all other paths.  This phenomenon, where decidability holds
for more complex sentences but not for simple ones, is part of
the inspiration for our concept of tree-decidability in Section \ref{sec:decidable}.

\begin{prop}
\label{prop:QE}
Let $T$ be a computable tree with computable finite branching and no isolated paths,
which (along with a functional $\Phi$) presents a torsion free abelian group $\A_T$.  Assume that every properly
existential-atomic formula defines a clopen subset of $\A_T$,
%that witnesses to such formulas are unique when they exist,
and that the decision procedure for this subset is uniform in the formula.
%and that we also have a procedure determining the unique solution
%(when one exists) uniformly in the formula and the parameters.
Then for every formula $\gamma(\Fvec)$ in the language of groups,
there exists a formula $\beta(\Fvec)$ in the same language
which is a Boolean combination of literals and subformulas
$\beta_i(\Fvec)$ such that every $\beta_i$ defines a clopen set and
$$ \A_T\models (\forall \Fvec)~[\gamma(\Fvec) \leftrightarrow \beta(\Fvec)].$$
Moreover, such a formula $\beta$, along with decision procedures
for each $\beta_i$, may be computed uniformly in $\gamma$.
\end{prop}
In Lemma \ref{lemma:Zplinear} we saw that the hypotheses here
apply to our presentation $T_p^+$ of $\Z_p^+$, for every prime $p$.
We will use it again in Section \ref{sec:decidable} for the group $\hatZ$.
We also clarify that this lemma does not claim full quantifier elimination,
as the subformulas $\beta_i$ could contain quantifiers.

\begin{pf}
We proceed by induction on the quantifier complexity of $\gamma(\Fvec)$,
which we take to be expressed in prenex normal form.
The proposition is trivial if $\gamma$ is quantifier-free.
If not, then suppose first that the proposition holds for up to $n$ quantifiers,
and suppose the $(n+1)$-st-innermost quantifier of $\gamma$
be $\exists G$.  (We handle the case $\forall G$ below.)
By induction we can re-express $\gamma$ to begin
with the same quantifier string up to $\exists G$,
followed by a matrix $\alpha(G,\Fvec,G_1,\ldots,G_k)$
that is a Boolean combination as described in the proposition.
Expressing this Boolean combination in disjunctive normal form,
we see that $\exists G\alpha$ is equivalent to an expression
$$ (\exists G\alpha_1)~\vee~\cdots~\vee~(\exists G\alpha_m),$$
where each $\alpha_j$ is a conjunction of literals and formulas $\beta_i(G,\Fvec,G_1,\ldots,G_k)$
defining clopen sets (for which we know decision procedures).
As $\A_T$ is abelian and torsion free, we may
perform the same reduction procedure as in Lemma \ref{lemma:reduced}
on each $(\exists G\alpha_j)$, moving any literals that do not involve $G$
outside the quantifier $\exists G$, but leaving the formulas $\beta_i$ inside it.

Within the scope of $\exists G$ after this operation,
there remains a conjunction of $\beta_i$'s and either a single equation
or finitely many inequations, all of which involve $G$ nontrivially.
If there are only inequations, then for every tuple $(\fvec,g_1,\ldots,g_k)$
from $\A_T$, densely many $g$ in $\A_T$ satisfy those inequations,
because otherwise $\A_T$ would have an isolated path.
In this case we preserve the conjunction of the $\beta_i$'s
(and the removable equations outside the quantifier $\exists G$)
but delete the quantifier $\exists G$ and the inequations, leaving an equivalent
expression that is exactly as required by the induction.
Each $\beta_i$ formerly defined a clopen set of elements $(g,\fvec,g_1,\ldots,g_k)$,
but since solutions $g$ to the inequations are dense, we now view $\beta_i$
as defining the projection of this clopen set onto $(\fvec,g_1,\ldots,g_k)$,
which is also clopen, uniformly in $\beta_i$.

On the other hand, if the scope of $\exists G$ contains a conjunction of $\beta_i$'s
and a single equation $\alpha^=$, then by the hypothesis of the proposition we may
determine a decision procedure for the (necessarily clopen) set $S$ defined
by $\exists G\alpha^=$.  Moreover, for each of the finitely many basic clopen sets
making up $S$, we may determine the corresponding initial segments of the
solutions $g$ to parameters in $S$.  By comparing these pairs with the
sets defined by the $\beta_i$'s, we may compute exactly the set of parameters
$(\fvec,g_1,\ldots,g_k)$ (without $g$) for which a solution $g$ of $\alpha^=$
exists that satisfies all the $\beta_i$'s.  This set of parameters is again clopen,
and our process yields a decision procedure for it, so the conjunction of 
$\exists G\alpha^=$ and the $\beta_i$'s is exactly what the proposition demands:
a formula defining a clopen set (of tuples $(\fvec,g_1,\ldots,g_k)$, which are the
remaining free and bound variables) with a decision procedure for that set.
Keeping the same literals along with this new formula completes the inductive step.

Finally, in case the $(n+1)$-st-innermost quantifier of $\gamma$ is $\forall G$
instead of $\exists G$, we simply express the inner formula $\forall G\alpha$
as $\neg\exists G(\neg\alpha)$, run the same procedure as above on $\exists G(\neg\alpha)$,
%(since the negation of a Boolean combination is still a Boolean combination
%and the negations of formulas defining clopen sets also define clopen sets),
and finally replace $\forall G\alpha$ by the negation of the new formula
given by that procedure, which is indeed of the desired form.

Thus, in both cases, we have eliminated the innermost quantifier from $\gamma$,
and have also eliminated the bound variable $G$.
Once this inductive step has been repeated enough times to eliminate all quantifiers
and all bound variables, the formula remaining instantiates the proposition.
\qed\end{pf}

Recalling Question \ref{question:Skolem}, we wish to determine which formulas have
computable (generalized) Skolem functions, given a particular tree presentation $T$.
Lemma \ref{lemma:reduced} yields a quick proof that in $T_p^+$, conjunctive existential formulas
all have computable generalized Skolem functions:  the reader will see how to compute
such a function for the reduced formula $\exists G\alpha'$ in the lemma, and this function
serves for the entire formula $\beta(\Fvec)\wedge\exists G\alpha'(\Fvec,G)$, since when
$\beta(\Fvec)$ is false the formula is unsatisfiable.  On the other hand, it is essential
that the original formula contain no disjunctions.  Indeed, from the following proposition,
we can see that there will always be purely existential formulas that have no computable Skolem functions. 

\begin{prop}
\label{prop:or issue}
Let $T$ be a tree presentation with at least one non-isolated path (in any signature whatsoever
with equality).  Let $\alpha(F,G)$ be the formula
$$\exists H[(F=G~ \wedge ~ F\neq H)~\vee~ (F\neq G ~\wedge ~ F=H) ]. $$
There is no computable Skolem function for $\alpha$ in $T$. 
\end{prop}

\begin{pf}
Let $P$ be a non-isolated path in $[T]$, and suppose that the Turing functional $\Psi$ computes a Skolem function for $\alpha$. We would then have $\Psi^{P,P}\neq P$. %but $\Psi^{P,Q}= P$ for any path $Q\neq P$. 
Take some $n$ such that $\Psi^{P,P}(n)\neq P(n)$. Since $\Psi^{P,P}$ must give an output after only viewing finitely much of its oracle, and since $P$ is non-isolated, there exists some $R\neq P$ such that $\Psi^{P,R}(n)=\Psi^{P,P}(n)\neq P(n)$, and so $\Psi^{P,R}$ fails to compute the unique witness $P$
to the (true) formula $\alpha(P,R)$.
\qed\end{pf}
Thus even purely existential formulas will not always have computable Skolem functions in computable tree presentations. However, if we limit our attention to formulas (of arbitrary quantifier complexity) that do not contain such a disjunction, then we can guarantee computable Skolem functions for a large class of computable tree presentations, including $T_p^+$. 

\begin{cor}
\label{cor:Skolem}
In the situation of Proposition \ref{prop:QE}, assume also
that witnesses to properly existential-atomic formulas are unique when they exist,
and that we also have a procedure determining the unique solution
(when one exists) uniformly in the formula and the parameters.  Then
every formula $\exists G\delta(F_1,\ldots,F_n,G)$ that can be expressed in prenex normal form
with a matrix that is a conjunction of literals (i.e., without the $\vee$ connective)
has a computable Skolem function $S$ for the tree presentation  $T$.
(That is:  given an oracle $\fvec$, $S$ computes a $g$ realizing
$\delta(\fvec,g)$ whenever one exists.)
\end{cor}
It should be noted that the hypothesis here applies to our tree presentation $T_p^+$ of $\Z_p^+$, for every prime $p$.
\begin{pf}
Apply Proposition \ref{prop:QE} to the formula $\exists G\delta(\Fvec,G)$.
A close examination of the procedure above shows that, after each step of
the induction, the formula remaining will still be in prenex form with a matrix
that is again a conjunction of literals and formulas defining clopen sets.
(The step for $\forall G$ uses $\neg\alpha$, which is a disjunction;
but it turns $\neg\alpha$ into a disjunction of literals and $\beta_i$'s
and then takes the negation of that, which is again a conjunction.)
At the end of the induction, we simply have a conjunction of such formulas.
The literals involve only the variables $\Fvec$, and it will in general
be undecidable whether $\fvec$ satisfies these literals.  However, we merely
check whether $\fvec$ satisfies all of the formulas defining clopen sets,
and if so, we use the uniqueness of solutions in the final step to produce
the desired $g$ (or else, if the final step had a conjunction of inequations in $F$,
produce any $g$ satisfying those inequations).  As long as $\fvec$ satisfies the
literals, this $g$ makes $\delta(\fvec,g)$ true.  If $\fvec$ fails to satisfy the
literals, then $\exists G\delta(\fvec,G)$ is false in the group, so by the definition
in Question \ref{question:Skolem}, it was OK to output any $g$ at all.
\qed\end{pf}

\section{Computing Isomorphisms}
\label{sec:iso}

In traditional computable structure theory, significant effort is devoted
to the computability of isomorphisms between computable structures.
For the purposes of the subject, two computable structures are considered
to be essentially identical when a computable isomorphism maps one
of them onto the other, since in this case the two computable structures share the
exact same computability-theoretic properties.  This contrasts with classical
model theory more generally, where any isomorphism at all between two structures
is usually deemed to make the two structures identical, even though their
computability-theoretic properties might be very different.  Our goal here is
to consider different computable tree presentations of the same structure $\A$
(that is, of two tree-presented structures identified with each other and
with $\A$ by classical isomorphisms).
In this case, each element of $\A$ can be named by a particular path
through a tree $T_0$ and also by a path through another tree $T_1$, with the
sets $[T_i]$ turned into copies of $\A$ by computable functions on these domains.
As always, a (classical) isomorphism between the presentations will simply mean a
bijection from $[T_0]$ onto $[T_1]$ respecting the functions in the signature.
The main question for us is when such an isomorphism should be called computable.

The simplest type of isomorphism arises from an isomorphism $f$ between
the trees themselves.  In this case, the isomorphism $f:T_0\to T_1$ naturally gives a map
$F:[T_0]\to[T_1]$.  Of course, this $F$ need not respect the functions in the signature,
but if it does, then the original $f$ may be seen as the germ of the isomorphism $F$.
In case $f$ is computable (as a map between the computable trees $T_0$ and $T_1$,
each of which is viewed as having domain $\omega$), then we would certainly
consider $F$ to be a computable isomorphism between the two presentations.
For an example of this situation, one need look no further than the two versions
of the group $\Z_p^+$ given before Definition \ref{defn:treepres}.  In the first,
a node in $T_0$ is just a tuple $(j_0,\ldots,j_{l-1})$ in $(\Z/p\Z)^l$ (for any $l\geq 0$),
ordered by the initial-segment relation.  In the second, a node in $T_1$ is an $l$-tuple
$(k_1,\ldots,k_l) \in (\Z/p\Z) \times (\Z/p^2\Z)\times\cdots\times (\Z/p^l\Z)$ (for any $l\geq 0$)
satisfying $(\forall n\leq l)(\forall m<n)~k_m\equiv k_n\bmod p^m$.  The natural isomorphism
between these trees maps $(j_0,\ldots,j_{l-1})$ to $(j_0,j_0+j_1p,\ldots,\sum_{i< l} j_ip^i)$,
and this is quickly seen to yield a bijection between $[T_0]$ and $[T_1]$ that is an isomorphism
from the group $([T_0],+)$ onto $([T_1],+)$, using the standard definition of addition
in each group.

The point of these sections is to transfer to $\Z_p\cross$
(and also to $\hatZ$ and $\prod_p\Z_p\cross$) the results in Section \ref{sec:existential}
about definable sets and Skolem functions for $\Z_p^+$.  The isomorphism we construct in 
Section \ref{sec:Zpcross} between $\Z_p\cross$ and $(\Z/(p-1)\Z)\times\Z_p^+$ will be of the nice sort
just described:  a computable isomorphism between the relevant trees that gives rise to a
computable isomorphism between the groups presented by those trees.
Consequently, our isomorphism will yield the desired results for $\Z_p\cross$,
once we complete Section \ref{sec:Zpcross}.  However, similar results hold in a broader
context, using the established notion of computable functions on Polish spaces, which we state here first.
\begin{defn}
\label{defn:iso}
Let $S$ and $T$ be trees presenting structures $\A_S$ and $\A_T$ in the same functional signature $\L$,
with computable operations $\Phi_F$ (in $S$) and $\Psi_F$ (in $T$) for each symbol $F$ in $\L$.
An isomorphism $H:\A_S\to\A_T$ is \emph{$C$-computable}
(for an oracle $C\subseteq\omega$) if there exists a Turing functional $\Theta$ such that
for every $f\in [S]$, $\Theta^{C\oplus f}$ is total and equals $H(f)$ in $[T]$.
Specifically, for each $n$-ary $F$ in $\L$ and each $n$-tuple $\fvec\in[S]^n$,
$$ \Psi_F^{\oplus_i (\Theta^{C\oplus f_i})} = 
\Theta^{C\oplus\Phi_F^{f_1\oplus\cdots\oplus f_n}}.$$
For a given $\Theta$ and $C$, we write $H^C_\Theta$ for the map from $[S]$ to $[T]$ defined thus.
\end{defn}

\begin{prop}
\label{prop:conjugation}
Let $S$ and $T$ be tree presentations of $\L$-structures.
Assume that there is a computable isomorphism $H$ from $\A_S$
onto $\A_T$ such that $H^{-1}$ is also computable.  If $\A_S$ has
a $C$-computable Skolem function for a given formula $\alpha$,
then $\B_S$ also has a $C$-computable Skolem function for that $\alpha$.
Moreover, the complexity of the subset of $(\A_S)^n$ defined by a formula
$\alpha(\Xvec,\Pvec)$ with parameters $\Pvec$ from $\A_S$ is exactly the complexity
of the subset of $(\A_T)^n$ defined there by $\alpha(\Xvec,H(\Pvec))$.
Finally, for any particular class $\mathfrak S$ of formulas, 
uniformity of these complexities and these Skolem functions
holds in $\A_S$ just if it holds in $\A_T$.
\end{prop}
\begin{pf}
This is clear.  The Skolem functions in $\A_T$ are just those in $\A_S$, conjugated by $H$.
\qed\end{pf}

The broad class of computable isomorphisms described by Definition \ref{defn:iso}
is not closed under inversion, so the requirement that $H^{-1}$ be computable
is nontrivial.  However, our isomorphism $F$ in Section \ref{sec:Zpcross}
will arise from a computable isomorphism $f$ between the trees themselves.
Therefore $f^{-1}$ will also be computable, and consequently so will $F^{-1}$.
(More generally, computable isomorphisms between structures presented by
trees with computable finite branching will always have computable inverses.)

\section{$(\Z_p)\cross$}
\label{sec:Zpcross}

It is natural to define a tree presentation of the multiplicative group $\Z_p\cross$
by taking the decidable subtree of the tree $T_p=(\Z/p\Z)^{<\omega}$ 
consisting of the root and all nodes whose first coordinate is not $0$.
The paths through this tree are precisely the elements of $\Z_p\cross$,
i.e., the units in $\Z_p$; the multiplication function from $T_p$ works
perfectly well on this subtree, and (multiplicative) inversion of paths
is computed coordinatewise modulo $p^n$ at each level $n$.
However, to prove results for $\Z_p\cross$ analogous to those
for $\Z_p^+$ in Section \ref{sec:existential}, it is easiest to use a different
tree presentation $T_p\cross$ of $\Z_p\cross$ that describes it using an isomorphism
$f$ (defined in Lemma \ref{lemma:range} below) involving the additive groups.

For all odd prime powers $p^n$, the multiplicative group $(\Z/p^n\Z)\cross$
is cyclic, containing those integers in $\{1,\ldots,p^n-1\}$
prime to $p$.  We define $v_n=p^n-p^{n-1}=p^{n-1}(p-1)$ to be its order.
In turn, this cyclic group is then isomorphic to the direct product
$(\Z/(p-1)\Z)^+\times(\Z/p^{n-1}\Z)^+$ of smaller cyclic groups (and $(\Z/(p-1)\Z)^+$
factors further, unless $p=3$).  Indeed, we will describe here a precise isomorphism
$\Z_p\cross \cong (\Z/(p-1)\Z)\times(\Z_p^+)$, which will enable us to present
the group $\Z_p\cross$ using of the results of Section \ref{sec:existential}.
(The case $p=2$ is a bit different, but our ultimate results will be similar.)

Assuming $p\neq 2$, we now facilitate consideration of the cyclicity of these groups.
%regard each multiplicative group $(\Z/p^n\Z)\cross$ as a cyclic group of order $v_n$.  To this end, 
Fix the least generator $q_1$ of the multiplicative group $(\Z/p\Z)\cross$.
For each $n\geq 1$, compute the least generator $q_{n+1}$ of the group $(\Z/p^{n+1}\Z)\cross$
with the property that 
$q_{n+1}\equiv q_n\bmod p^{n}$.  This is a finite uniform computation,
so we need only show that such a $q_{n+1}$ always exists.
Indeed, $q_n$ and all other generators lie in the same orbit under automorphisms
of $(\Z/p^n\Z)\cross$, and each such automorphism extends to an automorphism
of $(\Z/p^{n+1}\Z)\cross$, so each generator of $(\Z/p^n\Z)\cross$ is mod-$p^n$-equivalent
to exactly the same number of generators of $(\Z/p^{n+1}\Z)\cross$ as $q_n$ is.

\begin{lemma}
\label{lemma:extensions}
Let $x\in (\Z/p^n\Z)\cross$ and $y\in(\Z/p^{n+1}\Z)\cross$ with $y\equiv x\bmod p^n$.
Express these in terms of our generators as $x\equiv q_n^j\bmod p^n$ and $y\equiv q_{n+1}^k\bmod p^{n+1}$.
Then $j\equiv k \bmod v_n$.  (Recall that $v_n=p^n-p^{n-1}$.)
\end{lemma}
\begin{pf}
The generators satisfy $q_{n+1}^i\equiv q_n^i \bmod p^n$ for each power $i$.
Consequently, $q_{n+1}^i\equiv 1\bmod p^n$ just if $i$ is divisible by $v_n$ (the order of $q_n$).

We have $x=q_n^j \equiv q_{n+1}^j \bmod p^n$,
so $q_{n+1}^j \cdot q_{n+1}^{-k} \equiv xy^{-1} \equiv 1\bmod p^n$.
But %$q_{n+1}$ generates $(\Z/p^{n+1}\Z)\cross$, hence has order $v_n$.
with $q_{n+1}^{j-k}\equiv 1\bmod p^{n}$, our initial remark shows that
$v_n$ divides $(j-k)$.
\qed\end{pf}

\begin{cor}
\label{cor:extensions}
If $y_0=q_{n+1}^{k_0}$ and $y_1=q_{n+1}^{k_1}$ lie in $(\Z/p^{n+1}\Z)\cross$
with $y_0\equiv y_1\bmod p^n$, then $k_0\equiv k_1\bmod v_n$.
\end{cor}
\begin{pf}
Just apply Lemma \ref{lemma:extensions} once to each $y_i$ along with the unique $x\in(\Z/p^n\Z)\cross$
with $x\equiv y_0\bmod p^n$.
\qed\end{pf}

For all $n>0$, we define maps
$f_n : (\Z/p\Z)\cross\times (\Z/p^{n-1}\Z)^+ \to (\Z/p^n\Z)\cross$ by
$$ f_n( q_1^x, y) = q_n^{(xp^{n-1}+y(p-1))} \bmod p^n.$$
This is well-defined because $q_n$ has order $p^{n-1}(p-1)$ in $(\Z/p^n\Z)\cross$:
$$ f_n(q_1^{x+j(p-1)}, y+kp^{n-1})=q_n^{xp^{n-1}+j(p-1)p^{n-1} +y(p-1)+kp^{n-1}(p-1)}
\equiv f_n( q_1^x, y).$$
One checks readily that this $f_n$ is an isomorphism between the groups given.
The full map $f: (\Z/p\Z)\cross\times \Z_p^+ \to \Z_p\cross$ is then described by
$$ f(q_1^x,(y_1,y_2,\ldots)) = (f_1(q_1^x,0),f_2(q_1^x,y_1),f_3(q_1^x,y_2),f_4(q_1^x,y_3),\ldots).$$
This map $f$ respects the group operations.
We need to show that its image consists of actual $p$-adic integers.
\begin{lemma}
\label{lemma:range}
This map $f$ is a bijection from $((\Z/p\Z)\cross\times \Z_p^+)$ onto $\Z_p\cross$.
\end{lemma}
\begin{pf}
We first show that elements of the image satisfy the rule for belonging to $\Z_p$.
(They will then automatically be units there, as $f_1(q_1^x,0)=q_1^x\not\equiv 0\bmod p$.)
With $(y_1,y_2,\ldots)\in\Z_p^+$, we know that $p^{n-1}$ divides $(y_{n-1}- y_n)$
for every $n$.  Hence $(p-1)p^{n-1}$, which is the order $v_n$ of $q_n$, divides $((p-1)(y_{n-1}-y_n))$,
so
$$%q_{n+1}^{y_n(p-1)} \equiv 
q_n^{y_n(p-1)} \equiv q_n^{y_{n-1}(p-1)} \bmod p^n.$$
%with the first equivalence holding because $q_{n+1}\equiv q_n\bmod p^n$.
Turning to the other summand in the exponent, we see that
%know from the equivalence of the bases that $q_{n+1}^{xp^n} \equiv q_n^{xp^n}\bmod p^n$.  But now
$$ \left(q_n^{xp^n}\right)\cdot\left( q_n^{-xp^{n-1}}\right) = q_n^{x(p^n-p^{n-1})} \equiv 1\bmod p^n$$
because $q_n$ has order $(p^n-p^{n-1})$ in $(\Z/p^n\Z)\cross$.  Therefore
$q_n^{xp^n}\equiv q_n^{xp^{n-1}}\mod p^n$.  Combining all of this
and recalling $q_{n+1}\equiv q_n\bmod p^n$ yields
\begin{align*}
f_{n+1}(x,y_n) &\equiv q_{n+1}^{xp^n+y_n(p-1)}
\equiv q_{n}^{xp^n+y_n(p-1)} \bmod p^n\\
&\equiv q_{n}^{xp^n}\cdot q_{n}^{y_n(p-1)} \bmod p^n\\
&\equiv q_n^{xp^{n-1}} \cdot q_{n}^{y_{n-1}(p-1)} \equiv f_n(x,y_{n-1}) \bmod p^n.
\end{align*}
Thus $f$ does have image within $\Z_p\cross$.  Injectivity of $f$ is quickly proven
from its definition, noting the order of each $q_n$, and the surjectivity follows
because each $f_n$ is an injection between finite sets of equal cardinality.
\qed\end{pf}

It is important to notice that Lemma \ref{lemma:range} did not merely show $f$
to be an isomorphism of trees, but defined it so that it will give us a tree presentation of $\Z_p\cross$.
Recall from Section \ref{sec:existential} that we have computable tree presentations $T_p^+$
of each additive group $\Z_p^+$ of $p$-adic integers.  To produce our computable
tree presentation $T_p\cross$ of $\Z_p\cross$ (with $p>2$), we build the tree $T_p\cross$
which begins with $(p-1)$ nodes at level $1$, labeled from $1$ through $p-1$,
representing the elements of $(\Z/p\Z)\cross$, the first factor in the domain of $f$.
Each of these nodes at level $1$ in $T_p\cross$ then becomes the root of a copy
of the tree $T_p^+$, since $\Z_p^+$ is the second factor of that domain.
Formally, each node in $T_p\cross$ (except the root) is of the form $(x,\sigma)$,
with $x\in(\Z/(p-1)\Z)$ and $\sigma\in T_p^+$, and the partial order on these nodes
is given by
$$ (x,\sigma)\preceq (x',\sigma') \iff [ x=x'~\&~\sigma\sqsubseteq\sigma'].$$

We use the isomorphisms $f_n$ to label the nodes in $T_p\cross$.
Write $\lambda$ for the root of $T_p^+$. Each node $(x,\lambda)$ is labeled with
$q_1^x$, using our fixed generator $q_1$ for $(\Z/p\Z)\cross$.
Above that root, if $\sigma=(y_1,\ldots,y_{n-1})\in T_p^+$, then
the node $(x,\sigma)$ is labeled with $f_n(q_1^x,y_{n-1})$, an element of
$(\Z/p^n\Z)\cross$.  Thus, by Lemma \ref{lemma:range} and its proof,
the labels along each path in this tree $T_p\cross$ name an element of $\Z_p\cross$,
in such a way as to pair the paths bijectively with the units in $\Z_p$.
Under this identification, the group multiplication in $\Z_p\cross$ is
clearly computable:  to multiply two paths $g$ and $h$, one simply multiplies
(modulo $p^n$, for every level $n$ in the tree) the labels $g(n)$ and $h(n)$
of the nodes at level $n$ on these paths, and outputs the unique node at this level
in $T_p\cross$ with that label.  Likewise, the inverse of a path $g$
is computed just by taking the inverse in $(\Z/p^n\Z)\cross$ of each label $g(n)$.
So we have computable tree presentations $T_p\cross$ of each group
$\Z_p\cross$ (for $p>2$).

Next we consider the case $p=2$.  The groups $(\Z/2^n\Z)\cross$ are not cyclic,
except in case $n\leq 2$.  Rather, they are of the form 
$(\Z/2\Z)^+ \times (\Z/2^{n-2}\Z)^+$.  (This includes the case $n=2$, which
is cyclic of order $2$, while for $n=1$ the group itself is trivial.)
It is convenient that for $n>2$ we can use a fixed generator for the 
image in $(\Z/2^n\Z)\cross$ of the second factor.
\begin{lemma}
\label{lemma:5gen}
For every $n>2$, the element $5$ generates the subgroup
$G_n=\set{m\in(\Z/2^n\Z)\cross}{m\equiv 1\bmod 4}$, which is cyclic
of order $2^{n-2}$.
\end{lemma}
\begin{pf}
By inspection this holds for $n=3$, so assume inductively
that it holds for some $n\geq 3$.  The relevant subgroup $G_{n+1}$
of $(\Z/2^{n+1}\Z)\cross$ contains $2^{n-1}$ elements,
including $5$, so the order of $5$ there must be a power of $2$.
If $5$ were not a generator, then $5^{2^{n-2}}\equiv 1\bmod 2^{n+1}$,
so $2^{n+1}$ would divide $(5^{2^{n-2}}-1)$.  With $n\geq 3$, however,
$5^{2^{n-2}}-1 = (5^{2^{n-3}}-1)(5^{2^{n-3}}+1)$, and the latter factor
is $\equiv 2\bmod 4$, so $2^n$ must divide the first factor $(5^{2^{n-3}}-1)$.
Thus $5^{2^{n-3}}\equiv 1\bmod 2^n$, contradicting the inductive hypothesis.
\qed\end{pf}

Recognizing all this, we define isomorphisms
$f_n: (\Z/2\Z)^+ \times (\Z/2^{n-2}\Z)^+ \to (\Z/2^n\Z)\cross$ by 
$$ f_n(x,y) = (-1)^x\cdot 5^y \in \Z/2^n\Z$$
for all $n\geq 2$, while $f_1(x)=1$ maps $(\Z/2\Z)^+$ onto $(\Z/2\Z)\cross$.
We then put these together as before to create our desired isomorphism
$f: (\Z/2\Z)^+ \times \Z_2^+ \to \Z_2\cross$, defined (with $\yvec=(y_1,y_2,\ldots)\in\Z_2^+$) by
$$ f(x,\yvec) = (f_1(x),f_2(x,1),f_3(x,y_1),f_4(x,y_2)\ldots)=
(1,(-1)^x, (-1)^x\cdot 5^{y_1},\ldots)%(-1)^x\cdot 5^{y_2},\ldots)\in\Z_2\cross
.$$
It is clear that this map respects the group operations.  As before,
the main task is to show that it does produce elements of $\Z_2$
(which will clearly be units there, having initial coordinate $1$).
\begin{lemma}
\label{lemma:range2}
This map $f$ is a bijection from $((\Z/2\Z)^+ \times \Z_2^+)$ onto $\Z_2\cross$.
\end{lemma}
\begin{pf}
Assume that $\yvec\in\Z_2$, so $y_{n+1}\equiv y_n\bmod 2^n$ for every $n$.
Fix $n$ and express $y_{n-1}=y_{n-2}+k\cdot 2^{n-2}$.  Then for $n\geq 3$,
$$ f_{n+1}(x,y_{n-1})=(-1)^x 5^{y_{n-1}} = (-1)^x 5^{y_{n-2}}\cdot (5^{2^{n-2}})^k \equiv %5^{y_{n-2}}
f_n(x,y_{n-2})\bmod 2^n$$
by Lemma \ref{lemma:5gen}.  The cases
$n\leq 2$ are easily checked, so indeed $f(x,\yvec)\in\Z_2\cross$ as required.
As in Lemma \ref{lemma:range}, injectivity of the maps $f_n$ for $n>2$ is quickly seen;
$f_1$ itself is not injective, but the map $(f_1(x),f_2(x,y))$ maps the two
elements $(x,y)\in (\Z/2\Z)^+\times (\Z/2^{2-2}\Z)^+$ to the two distinct elements
$1$ and $3$ of $(\Z/2^2\Z)\cross$.  Surjectivity and bijectivity of $f$ itself then follow.
\qed\end{pf}

The tree $T_2\cross$ is now defined in the same general style as $T_p\cross$ was for odd $p$.
$T_2\cross$ has two nodes at level $1$, labeled by $0$ and $1$, the two elements
of the first factor $(\Z/2\Z)^+$ in the domain of $f$.  Each of these two nodes then
becomes the root of a copy of $T_2^+$, the tree presentation of the additive group
$\Z_2^+$, which is the second factor.  The multiplication and inversion in $T_2\cross$
are defined exactly as they were in $T_p\cross$ for odd $p$, using the isomorphism $f$
and the group operations on the factors.

The benefit of defining each $T_p\cross$ this way, rather than merely as a subtree of $T_p$,
is that it allows us to use answers to Questions \ref{question:elem}, \ref{question:definableset},
and \ref{question:Skolem} about $\Z_p^+$ to derive results about $\Z_p\cross$ as well,
most notably in Theorem \ref{thm:Zpcrosstreedecidable}.

It would of course be of serious interest to consider the ring $\Z_p$, as well as
its quotient field $\Q_p$, as opposed to 
considering $\Z_p^+$ and $\Z_p\cross$ separately.  It is well known,
from work of Ax and Kochen \cite{AK65} and independently Ershov \cite{E65},
that the theories of both this ring and its quotient field are decidable.
We leave these investigations for another time.
As remarked earlier, Kapoulas addressed our Question \ref{question:elem}
for $\Q_p$ in \cite{K99}.

\section{Cartesian products}
\label{sec:Cartesian}

To translate our results about $\Z_p^+$ into results in Section \ref{sec:decidable}
about $\Z_p\cross$, $\hatZ$, and $\prod_p\Z_p\cross$,
we need to relate facts about a direct product
of groups to facts about its factor groups.  We begin by explaining
tree presentations of direct products in general.

\begin{defn}
\label{defn:directproduct}
Let $T_0,\ldots,T_{n-1}$ be trees, each presenting a structure $\A_i=\A_{T_i}$
in the same language $\L$.  The \emph{direct product} of these tree presentations
is the subtree $T\subseteq\omega^{<\omega}$ containing those $\sigma$ such that,
for each $i<n$ and for $k =\lfloor\frac{|\sigma|-1-i}n\rfloor$,
the node $(\sigma(i),\sigma(i+n),\sigma(i+2n),\ldots,\sigma(i+kn))$ lies in $T_i$.
(The intuition is that we ``splice'' the trees $T_0,\ldots,T_{n-1}$ together,
level by level, to form $T$.  Level $nl+i$ of $T$ replicates level $l$ of $T_i$.)

For each $i$, if $\Psi_i$ represents the $r$-ary function symbol $f$ in $\A_i$, then
$f$ is represented in $\A_T$ by the functional $\Psi$ that accepts any oracle
$\oplus_{j\leq r}X_j$ of paths in $[T]$, decodes each $X_j$ into its sub-paths $X_{ji}\in[T_i]$,
runs $\Psi_i^{X_{1i}\oplus X_{2i}\oplus\cdots\oplus X_{ri}}$ to produce the path
$f^{\A_i}(X_{1i},\ldots,X_{ri})$, and splices these paths back together to form
the final output in $[T]$.  This simply says that $(X_{11},X_{12})\cdot (X_{21},X_{22}) =
(X_{11}\cdot X_{21},X_{12}\cdot X_{22})$, the standard operation on a Cartesian product.

For the Cartesian product of infinitely many trees $\la T_i\ra_{i\in\omega}$,
given uniformly in $i$, we do the exact same process, except that now
level $\la i,l\ra$ of $T$ replicates level $l$ of $T_i$, using a standard pairing function
$\la\cdot,\cdot\ra$ from $\omega^2$ onto $\omega$.
\end{defn}

In the specific case of $\Z_p\cross \cong (\Z/(p-1)\Z)^+\times \Z_p^+$,
the first factor is a finite group, and for simplicity we replicated the ``tree presentation''
of that factor simply by turning all its elements in nodes at level $1$, with a copy of $T_p^+$
(representing $\Z_p^+$) above each of those nodes.  This is effectively equivalent
to Definition \ref{defn:directproduct} in this case.  When we deal with
$\prod_{\text{primes~}p}\Z_p\cross$ and $\hatZ=\prod_{\text{primes~}p}\Z_p^+$,
we will use the definition's machinery for effective infinite Cartesian products.

Now it becomes necessary to relate facts about a direct product
of two groups to facts about its two factor groups.  We start with quantifier-free formulas.
In this section, when $G=\times_{i\in I}G_i$ (and $I$ could be infinite),
an $r$-tuple from $G$ consists of the rows of:
$$\xvec=\left(\begin{array}{ccccc}
x_{11} & x_{12} & \cdots & x_{1i} & \cdots\\
x_{21} & x_{22} & \cdots & x_{2i} & \cdots\\
\ldots & \ldots & \ddots & \ldots & \ddots\\
x_{r1} & x_{r2} & \cdots & x_{ri} & \ldots
\end{array}\right).$$
and we write $\xvec_i=(x_{1i},\ldots,x_{ri})$ for $r$-tuples in $(G_i)^r$;
and $\x_j = (x_{j1},x_{j2},\ldots)$ for single elements of the product $G$;
similarly for parameters.

\begin{lemma}
\label{lemma:CartesianQF}
Let $\alpha$ be a disjunction
$$\left(\vee_{j\leq m}\beta_j(\Xvec,\pvec)\right)\bigvee \left( \vee_{k\leq n} \neg\gamma_k(\Xvec,\pvec)\right),$$
in the language of groups, where all $\beta_j$ and $\gamma_k$
are atomic formulas with free variables $\Xvec$ and parameters $\pvec=(\pvec_0,\pvec_1)$
from a direct product $G=G_0\times G_1$ of two groups.
Fix a tuple $\xvec=(\xvec_0,\xvec_1)\in G^r$ of elements of $G$ as above.
Then $G\models \alpha(\vec x)$ if and only if
\begin{align*}
&G_0\models \vee_{k\leq n} \neg\gamma_k(\xvec_0,\pvec_0)
\text{~~or~~}
G_1\models \vee_{k\leq n} \neg\gamma_k(\xvec_1,\pvec_1)\\
&\text{~~or~~}\bigvee_{j\leq m} \Big[ G_0\models \beta_j(\xvec_0,\pvec_0)\text{~~and~~}G_1\models \beta_j(\xvec_1,\pvec_1)\Big].
\end{align*}
\end{lemma}
\begin{pf}
%Fix $\alpha$ and $\xvec$ as in the lemma.  
The important point is that
$$ G\models \beta_j(\xvec,\pvec)  \iff \left[G_0\models \beta_j(\xvec_0,\pvec_0)\text{~and~}
G_1\models \beta_j(\xvec_1,\pvec_1)\right],$$
whereas
$$ G\models \neg\gamma_k(\xvec,\pvec)  \iff \left[G_0\models \neg\gamma_k(\xvec_0,\pvec_0)\text{~or~}
G_1\models \neg\gamma_k(\xvec_1,\pvec_1)\right].$$
This makes it clear that the new expression is equivalent to $G\models\alpha$.
\qed\end{pf}

\begin{cor}
\label{cor:allQF}
For every quantifier-free formula $\alpha(\xvec,\pvec)$ and
every element $\xvec=(\xvec_0,\xvec_1)\in G=G_0\times G_1$,
there is a finite set of quantifier-free formulas $\delta_{0j}(\xvec_0,\pvec_0)$
and another finite set of such formulas $\delta_{1k}(\xvec_1,\pvec_1)$ such that
the truth of $\alpha(\xvec,\pvec)$ in $G$ is equivalent to
a Boolean combination of statements of the form $G_0\models\delta_{0j}(\xvec_0,\pvec_0)$
and statements of the form $G_1\models\delta_{1k}(\xvec_1,\pvec_1)$.
Moreover, this Boolean combination is computable uniformly in $\alpha$.
\end{cor}
\begin{pf}
Express $\alpha$ as a conjunction of formulas $\alpha_i$ in the form required by Lemma
\ref{lemma:CartesianQF}, i.e., disjunctions of literals, and apply that lemma
to get a conjunction of the Boolean combinations it gives for each $\alpha_i$.
%(Of course, the quantification ``for every $S\subseteq\{ 1,\ldots,m\}$'' can be converted
%into a conjunction of $2^m$ conjuncts, each with its own $S$.)
\qed\end{pf}

\begin{prop}
\label{prop:Cartesian}
For every formula $\alpha(\pvec)$ with parameters $\pvec=(\pvec_0,\pvec_1)$
from a direct product $G=G_0\times G_1$ of groups, the truth of $\alpha(\pvec)$ in $G$
is equivalent to a Boolean combination of statements $\Theta$ of the form
$G_i\models\delta(\Xvec_i,\pvec_i)$, which may be computed uniformly from $\alpha$.
For each $\Theta$, either $i=0$ or $i=1$,
so that each $\Theta$ considers only one of the factor groups, and
each $\delta$ is a formula in the language of groups, of no greater quantifier
complexity than $\alpha$.
\end{prop}
\begin{pf}
Corollary \ref{cor:allQF} provides the base case for an induction on the complexity of $\alpha$.
That case already completed the main task,
converting the quantifier-free portion of the statement $G\models\alpha(\pvec)$
into a Boolean combination of statements $\Theta$ of the desired form:
in particular, each $\Theta$ involves just one of $G_0$ and $G_1$.
It also showed how to compute the Boolean combination (and the statements
$\Theta$ themselves) uniformly from a quantifier-free $\alpha$.

For the inductive step, let $\alpha$ be $\exists X\beta(X,\pvec)$.  We use the chain of equivalences below,
where the inductive hypothesis states, for each $x=(x_0,x_1)\in G$, that $G\models\beta(x,\pvec)$ just if
$\vee_j\wedge_k\Theta_{jk}$ holds, with each $\Theta_{jk}$ as above.  (Here, for any fixed $x$, 
we treat $(x,\pvec)$ as the parameter tuple.  Fortunately the inductive procedure here is
uniform in the parameters.)  Notice also that for each $j$, $\wedge_k\Theta_{jk}$
is equivalent to a single conjunction
$$G_0\models \gamma_j(x_0,\pvec_0)~\&~G_1\models\delta_j(x_1,\pvec_1).$$
Indeed, two statements $G_0\models\phi(x_0,\pvec_0)$ and $G_0\models\psi(x_0,\pvec_0)$
both hold if and only if $G_0\models (\phi(x_0,\pvec_0)\wedge\psi(x_0,\pvec_0))$, and similarly for $G_1$,
so the conjunction $\wedge_k\Theta_{jk}$ can be shrunk down to the conjunction of just two such statements,
the first about $G_0$ and the second about $G_1$, exactly as written above.
Now all of the following hold.
\begin{align*}
G\models\exists X &\beta(X,\pvec) \iff (\exists x\in G)~G\models\beta(x,\pvec)\\
\iff& (\exists x_0\in G_0)(\exists x_1\in G_1)~G\models\beta((x_0,x_1),(\pvec_0,\pvec_1))\\
\iff& (\exists x_0\in G_0)(\exists x_1\in G_1)~\bigvee_j\bigwedge_k \Theta_{jk}\\
\iff& \bigvee_j(\exists x_0\in G_0)(\exists x_1\in G_1)~\bigwedge_k \Theta_{jk}\\
\iff& \bigvee_j(\exists x_0\in G_0)(\exists x_1\in G_1)~\left[(G_0\models\gamma_j(x_0,\pvec_0))~\&~ (G_1\models\delta_j(x_1,\pvec_1))\right]\\
\iff& \bigvee_j \left[(\exists x_0\in G_0)(G_0\models\gamma_j(x_0,\pvec_0))~\&~(\exists x_1\in G_1)
(G_1\models\delta_j(x_1,\pvec_1))\right]\\
\iff& \bigvee_j \left[(G_0\models\exists X_0~\gamma_j(X_0,\pvec_0))~\&~
(G_1\models\exists X_1~\delta_j(X_1,\pvec_1))\right].
\end{align*}
This last line is exactly as required by the Proposition.  Notice that, just as
$\alpha$ was more complex than $\beta$ by an initial $\exists$-quantifier,
so likewise the formulas in the last line are one $\exists$-quantifier more complex
than the $\gamma_j$ and $\delta_j$ corresponding to $\alpha$.

For $\alpha$ of the form
$\forall X\beta$, use the negations of the steps in this chain.
\qed\end{pf}

Corollary \ref{cor:allQF} and Proposition \ref{prop:Cartesian} both extend by induction
to finite Cartesian products of groups.  However, when we take a countable infinite
Cartesian product $G=\times_{i\in\omega} G_i$, the base case Lemma \ref{lemma:CartesianQF}
requires quantifiers.  It may be easiest to see the analogy using infinite $\vee$ and $\wedge$
connectives, as follows.  (The proof is just as in Lemma \ref{lemma:CartesianQF}.)

\begin{lemma}
\label{lemma:infCartesianQF}
Let $\alpha$ be a disjunction 
$$\left(\vee_{j\leq m}\beta_j(\xvec,\pvec)\right)\bigvee \left( \vee_{k\leq n} \neg\gamma_k(\xvec,\pvec)\right),$$
in the language of groups, where all $\beta_j$ and $\gamma_k$
are atomic formulas with free variables $\xvec$ and parameters $\pvec=(\pvec_i)_{i\in I}$
from a direct product $G=\times_{i\in I} G_i$ of infinitely many groups $G_i$.
Fix a tuple $\xvec\in G^r$ of elements
of $G$, and write $\xvec_i\in G_i^r$, for each $i$ such that $\xvec=(\xvec_i)_{i\in I}$.
Then $G\models \alpha$
if and only if
\begin{align*}
&\bigvee_{i\in I}~G_i\models \vee_{k\leq n} \neg\gamma_k(\xvec_i,\pvec_i)\\
\text{~~or~~}&\bigvee_{j\leq m} \bigwedge_{i\in I} G_i\models \beta_j(\xvec_i,\pvec_i).&\qed
\end{align*}
\end{lemma}
Equivalently, this says that $G\models \alpha$ if and only if
$(\exists i)~G_i\models \vee_{k\leq n} \neg\gamma_k(\xvec_i,\pvec_i)$
or
$\bigvee_{j\leq m} (\forall i)~G_i\models \beta_j(\xvec_i,\pvec_i).$
These appear to be more complex formulas, but really they are not.
After all, since each $G_i$ is tree-presented, the property
$G_i\models \neg\gamma_k(\xvec_i,\pvec_i)$
was already only $\Sigma^0_1$, while $G_i\models \beta_j(\xvec_i,\pvec_i)$
was already $\Pi^0_1$.  Adding the additional like quantifier to each
does not change the complexity.

Moreover, with $\alpha$ as above,
$\exists x_1\alpha(\xvec,\pvec)$ now becomes 
$(\exists i)~G_i\models \vee_{k\leq n} \exists x_1\neg\gamma_k(\xvec_i,\pvec_i)$
or
$\bigvee_{j\leq m} (\exists \vec{x_1}\in G)(\forall i)~G_i\models \beta_j(\xvec_i,\pvec_i)$,
and the second of these is equivalent to
$\bigvee_{j\leq m} (\forall i) (\exists x_{i,1}\in G_i)~G_i\models\beta_j(\xvec_i,\pvec_i)$,
hence also to
$\bigvee_{j\leq m} (\forall i)~G_i\models(\exists (x_1)_i)\beta_j(\xvec_i,\pvec_i)$.
So now we have a Boolean combination of statements $\Theta$ of the form
$G_i\models\delta(\Xvec_i,\pvec_i)$ where each $\delta$ is now existential,
i.e., of the same complexity as the formula $\exists x_1\alpha$ that we started with.
The same holds dually for $\forall x_1\alpha$.

\section{Tree-Decidability}
\label{sec:decidable}

In traditional computable structure theory, a countable structure with domain $\omega$
(in a computable signature) is said to be \emph{computable} if its atomic diagram
is decidable, and \emph{decidable} if its elementary diagram is decidable.  These
properties are completely natural to consider, as both diagrams are countable and
are readily coded into $\omega$ using G\"odel numberings.  The choice of
``computable'' for the one property and ``decidable'' for the other appears to have
been arbitrary, but it is well established and we see no reason to reconsider it.

For tree presentations, however, the situation is more complicated.
In a functional signature (with equality), we regard computable tree presentability
as a viable analogue, for continuum-sized structures, to computable presentability
for countable structures.  Indeed, we deem it to be the best analogue available,
unless one abandons the usual Turing model of digital computation
over a countable alphabet.
That said, the difficulty of handling relation symbols, even equality, has
long been apparent.  With uncountably many elements,
equality is inherently at least $\Pi^0_1$.  (To sketch the reason,
we expand on Lemma \ref{lemma:equality}:
a $\Sigma^0_1$ equivalence relation $\sim$ would decide $X\sim X$
using only finitely much information $\sigma$ from $X$, and then every $Y$
with the same information $\sigma$ would also have $Y\sim X$ under this
procedure.  Thus either some $Y\neq X$ has $Y\sim X$; or else
each $X$ is specified by some finite $\sigma$.  However, over a
countable alphabet there are only countably many finite strings $\sigma$ available.)

A second difficulty is the question of how to consider satisfaction of formulas
with quantifiers, as one does for countable structures by passing from computable
to decidable presentations.  For this issue, we apply a natural approach to
definable sets, formalized here as Definition \ref{defn:decidable}.
After presenting it, we will propose our preferred method of addressing
the problems related to equality.

\begin{defn}
\label{defn:decidable}
Let $T$ be a tree presentation of a first-order structure, and $\F$ a class of formulas
(potentially with free variables and first-order quantifiers) in its signature.  We say that
$T$ \emph{decides} the class $\F$ if there exists a Turing functional $\Phi$ such that,
for every formula $\phi(x_1,\ldots,x_n)$ from $\F$ (with all its free variables among $x_1,\ldots,x_n$)
and every $n$-tuple $(X_1,\ldots,X_n)\in [T]^n$ of elements of the structure,
the procedure
$$ \Phi^{X_1\oplus\cdots\oplus X_n}(\ulcorner \phi(\xvec)\urcorner)$$
halts and outputs $1$ if $\phi(X_1,\ldots,X_n)$ holds in the tree presentation $T$;
and halts and outputs $0$ if not.

This definition relativizes naturally to any fixed oracle set $C$,
defining \emph{$C$-decidability of $\F$ in the presentation $T$}, by using the procedure
$$ \Phi^{C\oplus X_1\oplus\cdots\oplus X_n}(\ulcorner \phi\urcorner).$$
\end{defn}

It is immediately seen that the singleton class $\{ x_1=x_2\}$ is not $C$-decidable this way
by any tree presentation of an uncountable structure, for any $C$.  This is a natural objection
to our definition.  On the other hand, the results in Section \ref{sec:existential} show
that the class of all nontrivial reduced existential formulas in the language of groups is decided by our tree
presentation $T_p^+$ of $\Z_p^+$, and also by our tree presentation of $\Z_p\cross$.
A similar result will be seen below for the field $\R$, when Definition \ref{defn:decidable}
is adapted to the situation of a tree quotient presentation.
Thus we consider the concept to be not without value.
It should be noted that both the theory of $\Z_p^+$ and the
theory of $\R$ are model-complete, with quantifier elimination down
to existential formulas.  We do not know at present whether any similar
results hold for structures without model-complete theories.

Definition \ref{defn:decidable} may be rephrased to say that every set
definable by a formula in $\F$ is a decidable subset of $[T]^n$, and that the
decision procedure is uniform in the formula.  This becomes somewhat
more impressive when one notices that formulas with parameters from $[T]$
fall under its umbrella, uniformly in the parameters,
provided there is no restriction on the number of free
variables in formulas in $\F$.  Indeeed, for a formula $\phi(x_1,\ldots,x_n,P_1,\ldots,P_m)$
with an $m$-tuple $\vec{P}$ of parameters from the presentation, one simply runs
$$ \Phi^{X_1\oplus\cdots\oplus X_n\oplus P_1\oplus\cdots\oplus P_m}(\ulcorner \phi(x_1,\ldots,x_{n+m})\urcorner)$$
on each $(X_1,\ldots,X_n)$ to decide the subset of $[T]^n$ it defines.

In light of this situation, we wish to offer Definition \ref{defn:treedecidable} below as
one new approach to the problem of deciding equality.  It simply recognizes the inherent
difficulty of the problem and incorporates that difficulty into the inputs to the decision procedure.

\begin{defn}
\label{defn:enumeration}
When $A\subseteq\omega$, an \emph{enumeration of $A$}
is simply a set $B\subseteq\omega$ such that $A=\pi_0(B)=\set{x}{\exists y~\la x,y\ra\in B}$,
where $\pi_0(\la x,y\ra)=x$ for all pairs $(x,y)\in\omega^2$,
using a standard computable bijection $(x,y)\mapsto\la x,y\ra$
from $\omega^2$ onto $\omega$.
In short, when $B$ is viewed as a subset of $\omega^2$, it should project onto $A$.

If such a $B$ is $\leq_T C$, we naturally call $B$ a $C$-computable enumeration of $A$
and say that $A$ is $C$-c.e.  This coincides precisely with the established meaning of the term.
\end{defn}
\begin{defn}
\label{defn:atomic}
Let $\L$ be a countable functional signature with equality, as always in this article,
and $\L'$ its extension by new constant symbols $c_1,c_2,\ldots$.
When $T$ is a computable tree presentation of an $\L$-structure
and  $(X_1,\ldots,X_n)\in[T]^n$, the \emph{positive atomic diagram
of $(X_1,\ldots,X_n)$ in $\A_T$} is the set
$$ \set{\ulcorner \alpha(c_1,\ldots,c_n)\urcorner}{\alpha\text{~atomic~\&~}
\A_T\models \alpha(X_1,\ldots,X_n)} \subseteq\omega,$$
containing those G\"odel codes of (positive) atomic formulas in $\L'$
realized in the substructure of $\A_T$ generated by $(X_1,\ldots,X_n)$.
\end{defn}
\begin{defn}
\label{defn:treedecidable}
Let $T$ be a computable tree presentation of a structure $\A$.
We say that the presentation $T$ is \emph{tree-decidable} if there exists a
Turing functional $\Phi$ that decides its elementary diagram
in the following sense:  whenever $\alpha(x_1,\ldots,x_n)$
is a formula from $\L$, and $(X_1,\ldots,X_n)\in[T]^n$
and $D$ is an enumeration of the positive atomic diagram
of $(X_1,\ldots,X_n)$ in $\A_T$,
$$ \Phi^{X_1\oplus\cdots\oplus X_n\oplus D}(\ulcorner\alpha(x_1,\ldots,x_n)\urcorner)\converges=
\left\{\begin{array}{cl}1,&\text{~if~}\A_T\models\alpha(X_1,\ldots,X_n);\\
0,&\text{~if~}\A_T\models\neg\alpha(X_1,\ldots,X_n).\end{array}\right.$$
\end{defn}

The intention of Definition \ref{defn:treedecidable} is to say that $T$ comes as close as
one can come to deciding (in the sense of Definition \ref{defn:decidable})
the elementary diagram of an uncountable structure.  The only snag
is the inevitable one of inability to decide equality of elements of the domain $[T]$.
We view tree-decidable structures as the best analogue, in power $2^\omega$,
to decidable structures in power $\omega$.  

It is common that ``most'' tuples
from $[T]^n$ will satisfy no nontrivial atomic formulas whatsoever.  (For example,
all $n$-tuples outside a set of measure $0$ in $\R^n$ are algebraically independent.)
A tuple that does satisfy such a formula is unusual.
The point of the definition is that, if our opponent demands
that we decide the elementary diagram of $\A_T$, then it is only fair that
the opponent should alert us when giving us a special tuple, such as
an $(X_1,X_2)$ in which $X_1$ happens to be the $37$-th power of $X_2$.
Indeed, we need only an enumeration of such unusual atomic facts
about the tuple, not a decision procedure, because we can enumerate
the negative atomic diagram ourselves, inequality being a $\Sigma^0_1$ relation
on $[T]$.  The real point is that, with the paths $X_1\oplus\cdots\oplus X_n$
and the positive atomic information that is given (plus the negative information
that we enumerate for ourselves using the paths), the procedure determines
the truth in $\A_T$ of arbitrarily complex first-order formulas, thus deciding
the elementary diagram of $\A_T$.

This definition admits obvious generalizations.  \emph{$C$-tree-decidability}
allows $\Phi$ to use an additional fixed oracle $C\subseteq\omega$; this might naturally
apply to a $C$-computable tree presentation of $\A$, for instance.  Restricting to a subclass
$\F$ of formulas requires $\Phi$ to succeed only on formulas from that class:
for example, one might speak of $\Sigma_n$ tree-decidability when $\F$ is the class
of all $\Sigma_n$ $\L$-formulas.  Finally, the definition generalizes to computable tree
quotient presentations, simply by requiring $\Phi$ to output $1$ if
the quotient structure $\A_T$ models $\alpha([X_1]_{\sim},\ldots,[X_n]_{\sim})$
and $0$ if not.

We close this article by showing that our ongoing examples (the tree presentations
$T_p$ of $\Z_p^+$ and $T_p\cross$ of $\Z_p\cross$) and certain products thereof
all satisfy Definition \ref{defn:treedecidable},
and that the generalization of Definition \ref{defn:treedecidable} to computable tree
quotient presentations is likewise satisfied by the standard example in that context,
namely the tree $T_{\R}$ of all fast-converging Cauchy sequences from $\Q$,
which presents the field $\R$ (in the signature of fields).

\begin{thm}
\label{thm:Zptreedecidable}
Our tree presentation $T_p^+$ of the group $\Z_p^+$, in the signature of additive groups with
equality, is tree-decidable.
\end{thm}
\begin{pf}
We give the procedure for $\Phi$, on a formula $\gamma(F_1,\ldots,F_n)$,
given parameters $\fvec\in(\Z_p^+)^n$ and the enumeration $D$ as described.
Proposition \ref{prop:QE} yields a formula, equivalent to $\gamma$, that is
a Boolean combination of literals and formulas $\beta_i$ defining known
clopen sets.  With the enumeration $D$ we can decide all the literals:
each equation will either be enumerated into $D$ eventually, or will
eventually be seen to be false at some level of the tree $T$; while
the reverse holds for each inequation.  Of course we can also determine
from the oracle for $\fvec$ whether each $\beta_i$ holds or not,
using the definitions of the clopen sets.  Thus we can decide the entire
Boolean combination, which is to say, we can decide the truth
of $\gamma(\fvec)$ in $\Z_p^+$.
\qed\end{pf}

\begin{thm}
\label{thm:Zpcrosstreedecidable}
The computable tree presentation $T_p\cross$ given in Section \ref{sec:Zpcross},
presenting the group $\Z_p\cross$, is tree-decidable.
\end{thm}
\begin{pf}
Recall that $T_p\cross$ codes $(\Z/(p-1)\Z)^+$ in its first level, followed by a copy of
$\Z_p^+$ above each level-$1$ node, corresponding to the direct product 
$(\Z/(p-1)\Z)^+\times\Z_p^+$, slightly modified when $p=2$.
Our procedure $\Phi$ for tree-decidability accepts a formula $\alpha(X_1,\ldots,X_n)$,
along with parameters $\cvec=((a_1,f_1),\ldots,(a_n,f_n))$ where all $a_i\in\Z/(p-1)\Z$
and $f_i\in\Z_p^+$, and uses Proposition \ref{prop:Cartesian} to convert the question
of the truth of $\alpha(\cvec)$ in $\Z_p\cross$ ($=\A_{T_p\cross}$) to that of a Boolean
combination of truth of certain formulas $\gamma(\avec)$ in $(\Z/(p-1)\Z)^+$
and truth of certain other formulas $\delta(\fvec)$ in $\Z_p^+$.  Since the former group
is finite, those questions can be answered trivially by $\Phi$, using $\avec$.

To decide whether $\Z_p^+\models \delta(\fvec)$ for a $\delta$ as above, it is natural
to appeal to Theorem \ref{thm:Zptreedecidable}, but there is one subtlety:
we need an enumeration of the atomic diagram $D$ of the subgroup of $\Z_p^+$
generated by the $\fvec$.  Of course we have an enumeration of the atomic diagram
$C$ of the subgroup of $\Z_p\cross$ generated by $\cvec$.  However, if (for example)
$c_1=(a_1,f_1)$ and $c_2=(a_2,f_1)$ with $a_2\neq a_1$, then $C$ will not
contain the (false) statement $c_1=c_2$, and so, to enumerate $D$, we must somehow
deduce that the second coordinates of these pairs are equal in $\Z_p^+$.  The following
lemma (with a small adjustment for the outlying case $p=2$) solves this difficulty.
\begin{lemma}
\label{lemma:enums}
In this context, for any two $(b_1,g_1),(b_2,g_2)\in (\Z/(p-1)\Z)^+\times\Z_p^+$, we have
(with $\Z_p\cross$ written multiplicatively)
$$ \Z_p^+\models g_1=g_2 \text{~if and only if~} \Z_p\cross\models (b_1,g_1)^{p-1}=(b_2,g_2)^{p-1}.$$
\end{lemma}
\begin{pf}
To see the forwards implication, notice that $(p-1)$ is the order of the first factor group.
Going backwards, we assume (writing
$\Z_p^+$ additively) that $(p-1)\cdot g_1 =(p-1)\cdot g_2$, and recall that, in the integral domain
$\Z_p$ of characteristic $0$, we have cancellation, making $\Z_p^+\models g_1=g_2$.
\qed\end{pf}
Using Lemma \ref{lemma:enums}, $\Phi$ can produce an enumeration of $D$ uniformly
from any enumeration of $C$ (enumerating not just those instances when $f_i=f_j$, but all
instances where a term in $\fvec$ equals $0$ in $\Z_p^+$).  This resolves the subtle point
and completes the proof of Theorem \ref{thm:Zpcrosstreedecidable}.
\qed\end{pf}

\begin{thm}
\label{thm:hatZtreedecidable}
The group $\hatZ$, the profinite completion of $\Z$, has a tree-decidable computable tree presentation.
\end{thm}
\begin{pf}
The presentation is simply the Cartesian product of the presentations $T_p^+$
of all the groups $\Z_p^+$, for all primes $p$, as in the version of
Definition \ref{defn:directproduct} for infinite products.  Notice that the subtlety
that arose in the proof of Theorem \ref{thm:Zpcrosstreedecidable} is no longer
so subtle:  our decision procedure for $\hatZ$, when given a tuple of parameters
$\bf{\vec c}=(\c_1,\ldots,\c_n)\in\hatZ^{<\omega}$ with each $\c_j=(c_{j2},c_{j3},c_{j5},\ldots)\in\prod_p\Z_p^+$,
is also given an enumeration of their atomic diagram in $\hatZ$, but two tuples
can be inequal even when many (or almost all) of their individual components are equal.
With the infinite product, we have no version of Lemma \ref{lemma:enums}
to which to resort.  Instead, we go through the same steps as for each individual $\Z_p^+$.

For a formula $\alpha$ of the form $\exists G~\sum a_iF_i=bG$, with $b\neq 0$, Corollary \ref{cor:cofinite}
shows that each group $\Z_p^+$ with $p$ not dividing $b$ contains a unique solution $g_p$
to $\exists G~\sum a_ic_{jp}=bG$, computable uniformly in $p$.
For the finitely many prime factors $p$ of $b$,
we apply Lemma \ref{lemma:Zplinear}:  if, for any of those primes, $\Z_p^+$
has no solution $G$ to the equation (with $(c_{1p},\ldots,c_{np})$ plugged in for
the $F_i$'s), then neither does $\hatZ$.  If all those primes
do yield solutions in $\Z_p^+$, then $\hatZ$ does have exactly one solution $G$,
the product of all the solutions in all the groups $\Z_p^+$, which we can compute
uniformly using the given parameters $\cvec$.  Thus the truth of simple existential formulas
of this sort in $\hatZ$ is decidable.  Moreover, for each such formula, the set it defines
is effectively clopen in $\hatZ$, uniformly in the formula.

For more complex formulas, we apply Proposition \ref{prop:QE}.
$\hatZ$ has now been seen to satisfy all the properties required there,
and the proposition gives effective quantifier elimination down to
Boolean combinations of literals and formulas defining clopen sets.
Just as with each group $\Z_p^+$ in Theorem \ref{thm:Zptreedecidable},
the enumeration $D$ of the atomic diagram of the subgroup generated
by the given parameters allows us to decide the truth in $\hatZ$ of the literals
for those parameters:  $D$ will tell us if equality holds, while inequality
is inherently $\Sigma^0_1$.  Also, the oracle for the parameters allows
us to decide their membership in the clopen sets resulting from
the quantifier elimination, and thus to decide whether the original formula
holds in $\hatZ$ with the given parameters.
\qed\end{pf}
It should be noted that this tree presentation of $\widehat{\Z}$ satisfies the hypothesis of Corollary \ref{cor:Skolem}, and so we get the computability of Skolem functions for formulas that are purely conjunctive (as described in the Corollary).

It is also good to see a situation in which tree-decidability fails.
For this we turn to the product of the multiplicative groups $\Z_p\cross$.
Understanding this theorem may make Theorem \ref{thm:hatZtreedecidable}
more impressive:  that proof rested on the fact that a nonzero coefficient $b$
can have only finitely many prime factors, while the analogous argument would fail here.
\begin{thm}
\label{thm:Zpcross}
The group $\mathcal C=\prod_p(\Z_p\cross)$, when presented as the direct product $T$
of the tree presentations $T_p\cross$ as in Definition \ref{defn:directproduct},
is not tree-decidable.  Indeed, even the class of existential formulas is not tree-decidable here.
\end{thm}
\begin{pf}
The basic point is that, as seen in Lemma \ref{lemma:range}, each factor
$\Z_p\cross$ has a subfactor $(\Z/p\Z)\cross$, which (for all $p>2$) has even order.
Therefore, the truth of the formula $\exists G~(F=G\cdot G)$ (written multiplicatively)
will fail to be uniformly decidable in $F$:  a parameter $f\in G$ that satisfies it
must have coordinates in every factor $(\Z/p\Z)\cross$ satisfying it, for every $p>2$,
but it can fail, in independent ways, in any of these factors $(\Z/p\Z)\cross$, leaving no way
to confirm that $\exists G~(f=G\cdot G)$ without checking all those coordinates.
(In odd-order groups, every element satisfies the formula.)

For a formal argument, let $f=(\vec 1,\vec 1, \vec 4,\vec 4,\vec 4,\ldots)\in\mathcal C=\prod_p\Z_p\cross$,
where $\vec 1$ represents the element $(1,1,1,\ldots)$ in each of $\Z_2$ and $\Z_3$,
while $\vec 4$ represents $(4,4,4,\ldots)$ in every other $\Z_p\cross$.
Clearly $\exists G~(f=G\cdot G)$ holds in $\mathcal C$, as witnessed by $(\vec 1,\vec 1, \vec 2,\vec 2,\vec 2,\ldots)$.
This $f$ has infinite order, so it generates an infinite cyclic subgroup of $\mathcal C$.
We let $D$ be a fixed computable enumeration of the atomic diagram of
the infinite cyclic group, and, given any Turing functional $\Phi$,
run $\Phi^{f\oplus D}(\ulcorner\exists G~(f=G\cdot G)\urcorner)$.
If $\Phi$ tree-decides $\mathcal C$, then this procedure must halt and
output $1$, having used only finitely much of its oracle $f\oplus D$.
Then we can run $\Phi$ on the same formula with a different oracle,
say $f'=(\vec 1,\vec 1, \vec 4,\vec 4,\vec 4,\ldots,\vec 4,x,\vec 4,\ldots)$,
and the same $D$, and get the same output $1$, provided that the $x$ here
lies in a factor $\Z_p\cross$ with $p$ so large that $x$ lies beyond the use of the
first computation.  When we do this with an $x=(y,z)\in\Z_p\cross=(\Z/p\Z)\cross\times(\Z_p^+)$,
choosing $y$ to be a generator of the cyclic group $(\Z/p\Z)\cross$, the same $D$
still enumerates the atomic diagram of the subgroup generated by $x$, but the formula
is now false in $\mathcal C$ for this $x$ (since the generator $y$ cannot be a square
in the even-order group $(\Z/p\Z)\cross$),
and so the output of $\Phi$ is incorrect in this case.
\qed\end{pf}

The preceding theorems also yield answers to Questions \ref{question:elem}
and \ref{question:definableset}.
Using Theorems \ref{thm:E72 1} and \ref{thm:E72 2}, one readily proves that
for every Turing ideal $I$, the groups $(\hatZ)_I$ and $(\prod_p \Z_p\cross)_I$
of $I$-computable elements are elementary subgroups of $\hatZ$ and
$\prod_p\Z_p\cross$ respectively.
For definable sets our answer to \ref{question:definableset} extends beyond
the context of groups.

\begin{prop}
\label{prop:definable}
When a computable tree presentation $T$ of a structure $\A_T$
is tree-decidable, every definable subset of $(\A_T)^n$ (for every $n$)
is arithmetically $\Delta^0_2$.
\end{prop}
\begin{pf}
If $\Phi$ is the procedure satisfying Definition \ref{defn:treedecidable},
then given any formula $\alpha(x_1,\ldots,x_n)$ and the jump $(\Xvec)'$
of any tuple $\Xvec\in(\A_T)^n$, we can run $\Phi$ on $\alpha$ using
the oracle $\Xvec$ (computed from $(\Xvec)'$) and deciding membership
of atomic sentences in $D$ when needed by using the full oracle $(\Xvec)'$.
\qed\end{pf}
It would be of interest to examine whether tree-decidability gives similar
straightforward answers to Questions \ref{question:elem} and \ref{question:Skolem}
in general.  One might also pose Question \ref{question:Skolem} in a context
where not just the parameters but also the relevant enumeration $D$
of their atomic diagram is given, as in Definition \ref{defn:treedecidable}.

The genesis of this project lay in questions about Galois groups of infinite algebraic
field extensions.  The two preceding theorems yield contrasting answers:
$\hatZ$, which is the absolute Galois group of every finite field, is tree-decidable
(thus, intuitively, very straightforward), while $\prod_p(\Z_p\cross)$, which is the
Galois group of the cyclotomic field over $\Q$, is more complex, despite being
abelian and relatively easy to describe.  For a further contrast, the article
\cite{M24} considers the absolute Galois group of $\Q$, which is viewed as a
far more complicated structure than either of these.  In contrast to $(\prod_p \Z_p\cross)_I$,
it is unknown whether the computable elements of $\operatorname{Gal}(\Q)$ form an
elementary subgroup.

As further examples, we now briefly demonstrate tree-decidability in two other known contexts.

\begin{thm}
\label{thm:Rtreedecidable}
The field $\R$, in the signature $(+,\cdot,0,1)$, has a 
tree-decidable tree quotient presentation.
\end{thm}
\begin{pf}
All the work is done by the Tarski-Seidenberg Theorem.  Given any
formula $\alpha(\xvec)$, that theorem effectively produces
a quantifier-free $\beta(\xvec)$, in the signature extended by $<$,
such that $\R\models\forall\xvec(\alpha\leftrightarrow\beta)$.  Next,
given the oracle $X_1\oplus\cdots\oplus X_n$ and the enumeration $D$
of the positive atomic diagram of the ring $\Z[\Xvec]$ (without $<$),
we proceed to decide the truth of each of the atomic
formulas used in $\beta$.  Each may be converted into either $0<h(\Xvec)$
or $0=h(\Xvec)$ for some $h$ in the polynomial ring $\Z[\Yvec]$.  In both cases,
we then run the two $\Sigma^0_1$ procedures that will halt just if
$0<h(\Xvec)$ and just if $h(\Xvec)<0$ (respectively), and simultaneously
watch for $D$ to enumerate the formula $0=h(\Xvec)$ as part of the
atomic diagram of the subfield $\Q(\Xvec)$ generated by these elements.
Exactly one of these three procedures will halt, and when we see that one halt,
we have decided the atomic formula in question.  Thus we may decide whether
$\R\models\beta(\Xvec)$, which is equivalent to $\R\models\alpha(\Xvec)$.
\qed\end{pf}

It may be noted that the procedure in Theorem \ref{thm:Rtreedecidable}
would succeed even in the signature $(+,\cdot,<,0,1)$ of ordered fields.
Additionally, we could have included subtraction and division among
the field operations here.  Doing all this would require an extension of
Definition \ref{defn:treepres}, which did not address either relation symbols
or partial functions (such as division).  However, given the enumeration $D$
of the atomic diagram in the original signature, we can decide the elementary diagram of $\R$
even in this extended signature.  It would be of interest to try to extend
Definitions \ref{defn:treepres} and/or \ref{defn:treedecidable}
to relational signatures more generally, using the concept of tree-decidability.
We will not attempt to do so here, but Theorem \ref{thm:Rtreedecidable}
suggests that the structure $\R$ in the signature of ordered fields, where
the relation $<$ is $\Sigma^0_1$, would be a natural example to consider.

We have not taken the space here to give a computable tree quotient presentation
of the field $\C$ of complex numbers:  the reader should be capable of doing this
without help, using $\R$ as a guideline.  Since $\C$ has full quantifier elimination,
effectively, without any use of additional relation symbols such as $<$, it will quickly
be seen that every reasonable computable tree quotient presentation
of $\C$ will also be tree-decidable.

Finally, we remind the reader that this article has deliberately given short shrift
to topological aspects of the situation.  We are computable structure theorists and wish
to advocate the consideration of tree-presented structures in computability.
However, when one moves beyond countable structures, topology becomes
an essential tool.  Definitions \ref{defn:treepres} and \ref{defn:treequotientpres}
clearly yield a computable topology on the structure $\A_T$ being presented,
and we expect that examining the interplay between model theory and that topology
will be far more fruitful in the long run than studying either one in isolation.
\\ \ \\
\textbf{Funding.}  Block and Miller were both supported by NSF grant \#DMS-2348891.
Additionally, Miller was supported by grant \#MP-TSM-00007933 from the
Simons Foundation and by several grants from the Research Foundation of CUNY.
%\\ \ \\
%\textbf{Acknowledgments.}


\begin{thebibliography}{99}

\bibitem{AK65}
J.\ Ax \& S.\ Kochen;
Diophantine problems over local fields II. A complete set of axioms for
$p$-adic number theory,
\emph{American Journal of Mathematics}
\textbf{87} 3 (1965), 631--648.

\bibitem{B24}
J.\ Block; Complexities of theories of profinite subgroups of $S_\omega$ via tree presentations, \emph{Twenty
Years of Theoretical and Practical Synergies} (2024),  366--376.

\bibitem{B25}
J.\ Block: Elementarity of subgroups and complexity of theories for profinite groups,
to appear in \emph{Computability}.

\bibitem{DH10}
R.G.\ Downey \& D.R.\ Hirschfeldt;
\emph{Algorithmic Randomness and Complexity}
(2010: Springer, New York).

\bibitem{DM25}
R.G.\ Downey \& A.\ Melnikov;
\emph{Computable Structure Theory:  A Unified Approach}
(2026: Springer, New York).

\bibitem{DM23}
R.G.\ Downey \& A.\ Melnikov; Computably compact metric spaces, \emph{The Bulletin of Symbolic
Logic} \textbf{29} 2 (2023), 170-–263. 

\bibitem{E72}
P.C.\ Eklof;
Some model theory of Abelian groups, 
\emph{Journal of Symbolic Logic}
\textbf{37} 2 (1972), 335--342.

\bibitem{E65}
Y.L.\ Ershov;
On the elementary theory of maximal normed fields,
\emph{Doklad.\ Ak.\ Nauk.\ USSR} 
\textbf{6} (1965),1390--1393.

\bibitem{FJ86}
M.D.\ Fried \& M.\ Jarden;
\emph{Field Arithmetic}
(Berlin:  Springer-Verlag, 1986).

\bibitem{H23}
B.\ Hart;
An introduction to continuous model theory;
chapter in \emph{Model Theory of Operator Algebras},
ed.\ I.\ Goldbring (Berlin:  De Gruyter, 2023).

\bibitem{K99}
G.\ Kapoulas;
Computable $p$-adic numbers,
CDMTCS Research Report Series, 1999.

\bibitem{K95}
A.S.\ Kechris;
\emph{Classical Descriptive Set Theory}
(New York:  Springer-Verlag, 1995).


\bibitem{K23}
H.T.\ Koh, A.\ Melnikov, \& K.M.\ Ng;
Computable topological groups, 
\emph{The Journal of Symbolic Logic}
\textbf{90} (2025) 1, 188--220.

\bibitem{KK17}
M.\ Korovina \& O.\ Kudinov; 
Spectrum of the computable real numbers,
\emph{Algebra and Logic}
\textbf{55} (2017) 6, 485--500.

\bibitem{L81}
P.\ La Roche;
Effective Galois theory,
\emph{Journal of Symbolic Logic}
\textbf{46} 2 (1981), 385--392.


\bibitem{M24}
R.\ Miller;
Computability for the absolute Galois group of $\Q$,
submitted for publication.

\bibitem{P20}
A.\ Pauly, D.\ Seon, \& M.\ Ziegler: Computing Haar measure, \emph{28th EACSL Annual
Conference on Computer Science Logic} 
 (2020), 34:1--34:17.


\bibitem{S81}
R.L.\ Smith;
Effective aspects of profinite groups,
\emph{Journal of Symbolic Logic}
\textbf{46} 4 (1981), 851--863.

\bibitem{S87}
R.I.\ Soare;
\emph{Recursively Enumerable Sets and Degrees}
(New York:  Springer, 1987).


\bibitem{S55}
W.\ Szmielew;
Elementary properties of Abelian groups
\emph{Fundamenta Mathematicae} \textbf{41} (1955), 203--271
DOI: 10.4064/fm-41-2-203-271.

\bibitem{W00}K.\ Weihrauch; \emph{Computable Analysis} (Berlin:, Springer, 2000).

\bibitem{W09}
K.\ Weihrauch \& T.\ Grubba;
{Elementary Computable Topology}
\emph{Journal of Universal Computer Science} \textbf{15} 6 (2009), 1381--1422.

\end{thebibliography}
\end{document}